\documentclass{amsart}
\usepackage{amsfonts}

\usepackage{graphicx}
\usepackage{psfrag}

\newtheorem{theorem}{Theorem}[section]
\theoremstyle{plain}
\newtheorem{corollary}[theorem]{Corollary}
\newtheorem{lemma}[theorem]{Lemma}
\newtheorem{proposition}[theorem]{Proposition}
\theoremstyle{remark}
\newtheorem{remark}[theorem]{Remark}

\numberwithin{equation}{section}

\newcommand{\tr}{\operatorname{tr}}
\newcommand{\otr}{\operatorname{0-tr}}
\newcommand{\ovol}{\operatorname{0-vol}}
\newcommand{\FP}{\operatornamewithlimits{FP}}

\newcommand{\re}{\operatorname{Re}}

\newcommand{\res}{\operatorname{Res}}
\newcommand{\rank}{\operatorname{rank}}
\newcommand{\ord}{\operatorname{ord}}

\newcommand{\bbR}{\mathbb{R}}
\newcommand{\bbH}{\mathbb{H}}
\newcommand{\bbC}{\mathbb{C}}
\newcommand{\bbZ}{\mathbb{Z}}
\newcommand{\bbN}{\mathbb{N}}

\newcommand{\cinf}{C^\infty}

\newcommand{\bX}{{\partial_\infty X}}
\newcommand{\FPe}{\FP_{\varepsilon\downarrow 0}}
\newcommand{\Oint}{\sideset{^{0\hskip-5pt}}{}\int}

\begin{document}

\title[Selberg's zeta function and spectral geometry]{Selberg's zeta function and the 
spectral geometry of geometrically finite hyperbolic surfaces}
\author[Borthwick]{David Borthwick}
\address[Borthwick]{ Department of Mathematics and Computer Science, Emory
University, Atlanta, Georgia, 30322, U. S. A.}
\email{davidb@mathcs.emory.edu}
\author[Judge]{Chris Judge}
\address[Judge]{ Department of Mathematics, Indiana University, Bloomington,
Indiana, 47401. U. S. A.}
\email{cjudge@indiana.edu}
\author[Perry]{Peter A. Perry}
\address[Perry]{ Department of Mathematics, University of Kentucky,
Lexington, Kentucky, 40506--0027, U. S. A.}
\email{perry@ms.uky.edu }
\thanks{Borthwick supported in part by NSF\ grant DMS-0204985.}
\thanks{Perry supported in part by NSF\ grant DMS-0100829.}
\date{November, 2003}

\begin{abstract}
For hyperbolic Riemann surfaces of finite geometry, 
we study Selberg's zeta function and its relation to
the relative scattering phase and the resonances of the Laplacian.  
As an application we show that the conjugacy class of a finitely 
generated, torsion-free,  discrete subgroup of $SL(2, {\mathbb R})$ 
is determined  by its trace spectrum up to finitely many possibilities, 
thus generalizing results of McKean \cite{McKean:1972} and
M\"uller \cite{Mueller:1992} to groups 
which are not necessarily cofinite.
\end{abstract}

\maketitle

\section{Introduction}\label{intro.sec}

A geometrically finite hyperbolic surface is a topologically finite, 
complete Riemannian surface of constant curvature -1. 
In this paper we study the Selberg zeta function, $Z_X(s)$
(see \S \ref{zeta.sec} for the definition),
associated to the length spectrum of a geometrically finite surface $X$.   
When $X$ has infinite area, the discrete spectrum of the
positive Laplacian $\Delta _{X}$ will be finite and possibly empty.
The appropriate spectral invariants at this level of generality are the 
{\em scattering resonances} or simply {\em resonances}.  These are the
poles of the meromorphic continuation of the resolvent 
$(\Delta_X - s(1-s))^{-1}$,
with multiplicities (see \S\ref{scatt.sec} for the definition).

Building upon our previous results  and the work of Guillop\'{e} and Zworski 
\cite{GZ:1997},  we exhibit an explicit connection (Theorem 
\ref{thm.zeta2}) between Selberg's zeta 
function and a Hadamard product over the resonance set.  
This yields the following:

\begin{theorem}
\label{thm.zeta} 
Let $n_C$ denote the number of cusps of $X$,
and let $\chi$ denote the Euler characteristic of $X$. 
The Selberg zeta function $Z_{X}$ extends to a meromorphic function of
order two, with a divisor that can be divided into spectral and topological
components:  The spectral zeros of $Z_X$ 
are given by the resonance set $\mathcal{R}_X$ (with multiplicities).
In addition, $Z_X(s)$ has topological zeros at $s = -k$ for $k \in \bbN_0$,  
of order $(2k+1) \cdot (-\chi)$, and topological poles of order $n_C$ 
at $s \in\frac12 - \bbN_0$.
\end{theorem}

In the full geometrically finite context, the meromorphic continuation of 
$Z_{X}$ to $\mathbb{C}$ was proven by Guillop\'{e} \cite{Guillope:1992}. 
If $X$ has finite area, then meromorphic continuation with the order 
bound and divisor as given here can be deduced from the Selberg trace formula
(see \cite{He:1983}).
In the convex co-compact case ($n_{C}=0$) Patterson-Perry  \cite{PP:2000}
established the order bound, using thermodynamic formalism  
of Ruelle and Fried, and also computed the divisor.
The thermodynamic formalism does not extend to the case $n_C >0$.
To prove that Theorem \ref{thm.zeta} holds for all geometrically finite $X$ 
requires a different approach from these previous results.

As in the compact case, one can exploit knowledge of the divisor of the zeta
function to link the resonance set to the length spectrum.

\begin{corollary}\label{cor.finite.euler}
The resonance set determines the length spectrum of $X$, 
the Euler characteristic, $\chi$, and the number of cusps, $n_C$.  
The length spectrum determines $\chi$ and $n_C$ up to a finite
number of possibilities.  The length spectrum, $\chi$, and $n_C$,
together determine the resonance set.
\end{corollary}
At present, we do not know whether  the
length spectrum (by way of $Z_X$) determines $\chi$ and $n_C$ when  
$X$ has infinite area and $n_C>0$.  That is, we cannot rule out the 
possibility that the resonance set contributes
to the multiplicity of each of the zeros of $Z_X$ that lie in $1/2-{\mathbb N}_0/2$.
However, if $X$ has finite area, the resonances are confined to 
a vertical strip, and hence for sufficiently negative $k \in -{\mathbb N}_0$,
the multiplicity of the zero of $Z_X$ equals $(2k+1) \cdot (-\chi)$.  
And in the convex co-compact
case (infinite area with no cusps), Martin Olbrich has pointed out to us that 
Corollary 6.9 of \cite{Olbrich:2002}, in conjunction with the theory of 
\cite{BU:1999},
shows that the multiplicity of resonances at $-\bbN$ is either
zero (for non-elementary groups) or two (for elementary groups).  Thus,
in the convex co-compact case $\chi$ can be recovered as
$(\ord_{s=-k} Z_X - \ord_{s=-k-1} Z_X)/2$ for any $k\in\bbN$.
\footnote{The possibility 
of resonances overlapping topological zeros was was overlooked in 
\cite{BJP:2002}, but Olbrich's result implies that Theorem 3.1 of that paper is
correct as stated.}

By combining  Corollary \ref{cor.finite.euler} with methods of Teichm\"uller
theory we obtain the following application:

\begin{theorem}\label{Finiteness}
Let $X$ be a complete, geometrically finite hyperbolic surface of 
infinite area.  Then the length spectrum of $X$ 
determines $X$ up to finitely many possibilities.
In particular, the resonance set determines $X$
up to finitely many possibilities.
\end{theorem}
Our proof of the first claim in Theorem \ref{Finiteness}
requires Theorem \ref{thm.zeta}. Indeed, we do not know how to obtain
bounds on the Euler characteristic without using scattering theory.

If $X$ is connected, then there exists a finitely generated, torsion free, 
discrete subgroup, $\Gamma$, of $SL(2, {\mathbb R})$, unique up to conjugation,
such that  $X$ is the quotient of the hyperbolic upper half-plane, $\bbH$,
by $\Gamma$ acting as M\"obius transformations. 
Conversely, any such quotient is  a geometrically finite hyperbolic surface.  
Moreover,  the length spectrum of  $X$ equals (twice) the set 
$\{ \cosh( \rm{trace}(\gamma))~ |~ \gamma \in \Gamma\}$ 
(including multiplicities). Thus, from Theorem \ref{Finiteness} we have 

\begin{corollary} \label{Discrete}
Let $\Gamma$ be a finitely generated, torsion-free, 
discrete subgroup of  $SL(2, {\mathbb R})$ of co-infinite area.  Then
the set of traces $\{ \rm{trace}(\gamma)|~ \gamma \in \Gamma\}$
(including multiplicities) determines the conjugacy class
of $\Gamma$ up to finitely many possibilities.
\end{corollary} 

Theorem \ref{Finiteness} is due to McKean \cite{McKean:1972}
in the case that $X$ is compact and to M\"uller \cite{Mueller:1992} in the
case that $X$ has finite area.   Corollary \ref{Discrete} is also due
to McKean \cite{McKean:1972} in the case that $\Gamma$ is
cocompact.   

\begin{figure} \label{gfsurface}
\caption{Decomposition of the surface $X$}

\psfrag{C1}{$C_1$}
\psfrag{C2}{$C_2$}
\psfrag{Z}{$Z$}
\psfrag{F1}{$F_1$}
\psfrag{bX}{$\bX$}
\psfrag{horocycle}{horocycle}
\psfrag{geodesic}{geodesic}
\begin{center}   \includegraphics{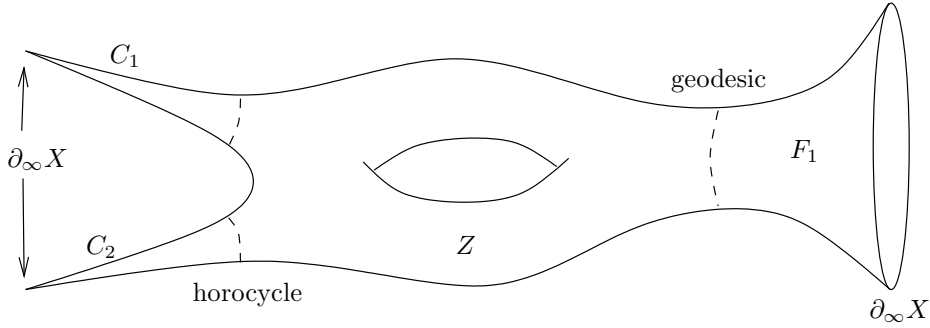} 
\end{center}
\end{figure}

Theorem \ref{thm.zeta} is proven by linking the zeta function to 
the scattering theory of $\Delta_X$.  For 
geometrically finite $X$ scattering theory can be set up
(following the approach of Guillop\'{e} and Zworski 
\cite{GZ:1997, GZ:1999}) using a decomposition of the form 
\begin{equation}\label{zcf.decomp}
   X~ =~ Z\sqcup \left( C_{1}\sqcup \cdots \sqcup C_{n_{C}}\right) \sqcup 
   \left(F_{1}\sqcup \cdots \sqcup F_{n_{F}}\right)
\end{equation}
where $Z$ is a compact surface,  each  $F_i$, $1 \leq i \leq  n_F$,
is a funnel, and each $C_j$, $1 \leq j\leq n_C$ is a cusp (end).
This is illustrated in Figure \ref{gfsurface}.
(For an exact description of funnels and cusps, see \S  \ref{sec.model}.)
The boundary of $Z$ consists of $n_{F}$ closed geodesics 
(uniquely determined) and $n_{C}$ 
horocycles (the choice of which is not unique) along
which $Z$ is glued to the funnel and cusp ends, respectively. 
The ideal boundary $\bX$ is the disjoint union of $n_C$
ideal points and $n_F$ circles.  

The scattering operator $S_{X}(s)$, described in greater detail in
\S\ref{scatt.sec}, can be viewed as a map from 
$C^\infty(\bX)$ to itself, where a ``smooth function'' on $\bX$ 
is interpreted as an element of $\mathbb{C}^{n_{C}}\oplus \left(
\oplus _{j=1}^{n_{F}}C^{\infty }(S^{1})\right)$. 
Because the funnels are isometric to standard models (see \S \ref{sec.model}),
one can define relative scattering from the disjoint union
\begin{equation*}
Y=F_{1}\cup F_{2}\cup \cdots \cup F_{n_{F}}
\end{equation*}
to $X$. On each funnel end, the Laplacian $\Delta _{F_{j}}$ carries
Dirichlet boundary conditions on the closed geodesic which joins $F_{j}$ to 
$Z$.   The scattering operator for $\oplus _{j=1}^{n_{F}}\Delta _{F_{j}}$, denoted 
$S_{Y}(s)$, is a direct sum of scattering operators for the
Dirichlet Laplacians $\Delta _{F_{j}}$, and acts on $\oplus _{j=1}^{n_{F}}
C^{\infty }(S^{1})$. If we set
\begin{equation*}
S_{0}(s)=\left( \oplus _{j=1}^{n_{C}}\mathbf{1}\right) \oplus 
S_{Y}(s)
\end{equation*}
(where $\mathbf{1}$ is the identity map on $\mathbb{C}$) then $S
_{0}(s)$ acts on $C^{\infty }(\bX)$ and the relative scattering operator
\begin{equation*}
S_{rel}(s)=S_{X}(s)S_{0}(s)^{-1}
\end{equation*}
is determinant class.  Thus one has a relative scattering determinant: 
\begin{equation*}
\tau _{X}(s)=\det \left( S_{rel}(s)\right) .
\end{equation*}

In \cite{GZ:1997}, Guillop\'e and Zworski show that $\tau_{X}(s)$ extends to the 
quotient of analytic functions of order at most four and compute its divisor in
terms of the resonance set.  
Our proof of Theorem \ref{thm.zeta} in \S\ref{zeta.sec} involves establishing
the connection between $Z_X(s)$ and $\tau_X(s)$, and then exploiting this
information about $\tau_X(s)$.
In the course of the proof we will improve the order bound on $\tau_X(s)$
from four to two.

\textbf{Acknowledgments}. We are grateful to Reza Chamanara and 
 Jeffrey McGowan for helpful conversations, and to Martin Olbrich for point out
to us his results in \cite{Olbrich:2002}. We also thank Peter Sarnak for 
encouraging us  to remove the unnecessary geometric 
assumptions used in \cite{BJP:2002}.



\section{Regularized traces and model problems}\label{sec.model}

We start with $X$ a geometrically finite hyperbolic surface 
decomposed into a compact core $Z$ plus cusps $C_i$ 
and funnels $F_j$ as in (\ref{zcf.decomp}). 
The cusp and funnel ends are isometric to
standard models which we now describe.
Here and in what follows, we use the upper half-plane model 
of hyperbolic space: $\bbH = \{z=x+iy| \>y>0\}$, with metric
\begin{equation*}
          ds^{2}=y^{-2}(dx^{2}+dy^{2})\text{.}
\end{equation*}

For $\ell>0$, denote by $\Gamma_{\ell}$ the abelian 
discrete group generated by the hyperbolic isometry $z \mapsto e^{\ell} \cdot z$.  
The quotient  $\Gamma_\ell \backslash  \bbH$ is the 
{\em hyperbolic cylinder} with a single primitive closed geodesic of 
length $\ell$.  It is well-known and easy to see that this cylinder is isometric 
to ${\mathbb R} \times  {\mathbb R}/{\mathbb Z}$ equipped with 
the metric
\begin{equation} 
        dr^{2}+\ell^2 \cdot
           \cosh ^2(r)~  dx ^{2},  \label{eq.Y.metric}
\end{equation}
A {\em funnel} is the subspace $[0, \infty) \times {\mathbb R}/{\mathbb Z}$ 
equipped with this metric. In particular, the boundary, $\partial F$,  of $F$
is a geodesic of length $\ell$.

Let $\Gamma _{\infty }$ the abelian discrete group generated by the 
parabolic isometry  $z\mapsto z+1$.  For us  the  {\em horn} is the 
quotient 
\begin{equation*}
       H~ =~ \Gamma _{\infty }\backslash {\bbH}.
\end{equation*}
It is well-known and easy to see that the horn is isometric to
${\mathbb R} \times {\mathbb R}/{\mathbb Z}$ equipped
with the metric 
$$
             dr^2~  +~  e^{-2r}~ dx^2.
$$
A {\em cusp end} is a subspace of the form $[a, \infty) \times 
{\mathbb R}/{\mathbb Z}$ equipped with this metric.
A level curve of the coordinate $r$ is called a {\em horocycle}.
In particular, a cusp end is bounded by a horocycle of $r=a$.
Note that in the decomposition of (\ref{zcf.decomp})
the parameter $a$ can be chosen to be $0$ 
independent of the surface. This follows from a 
collar type lemma for cusps to be found, for example, in \cite{Buser:1992}.

We'll regularize integrals over $X$ by
introducing a `defining function' $\rho$ for the boundary at infinity,
Namely, choose $\rho  \in  \cinf(X)$ so that in each of the cusp 
and funnel ends $\rho =e^{-r}$ with respect to the geodesic coordinate
systems introduced above.   With this choice we can define a regularized integral:   
Suppose $f\in \cinf(X)$ is \textit{polyhomogeneous} in $\rho$, meaning that it
has an asymptotic expansion as $\rho\to0$ in powers of $\rho$ and $\log\rho$,
with smooth coefficients.  We then define 
\begin{equation*}
\Oint_{X}f\>dg=\FPe
\int_{\{\rho \geq \varepsilon \}}f\>dg,
\end{equation*}
where $FP$ denotes the finite part, meaning that the limit is taken after the 
divergent terms in the asymptotic expansion as $\varepsilon\to 0$ are 
subtracted off.  The 0-volume of a region is correspondingly defined as the 
0-integral 
of $dg$.  These regularizations would depend on the choice of $\rho$ in general. 
Our definition exploits the canonical choice of $\rho$ made possible for 
hyperbolic surfaces  by the standard models for the ends.

\begin{lemma}\label{ovol.lemma}
With the $0$-integral as defined above,
\begin{equation*}
\ovol(X)~ =~ -2\pi \cdot \chi,
\end{equation*}
where $\chi$ is the Euler characteristic of $X$
\end{lemma}
\begin{proof}
The area of the region $\{\rho \geq \varepsilon \} \subset F$ is 
$\ell (\varepsilon ^{-1}-\varepsilon )/2$, so $\ovol(F)=0$.  Hence $\ovol(X)
= \operatorname{vol}(X\setminus Y)$ and the result follows from 
the Gauss-Bonnet theorem, since the funnels have geodesic boundary.
\end{proof}

For a smoothing operator $T$ defined on $X\times X$, with continuous 
kernel $K_{T}(x,y)$ with respect to Riemannian measure on $X$, 
the $0$-trace of $T$ is defined by
\begin{equation*}
\otr T = \Oint_{X} K_{T}(z,z)\,dg(z).
\end{equation*}
For this to be well-defined, the function $K_T(x,x)$ must be polyhomogeneous
as a function of $\rho$, which will be the case for the operators we consider.

There are two regularized traces formed from the resolvent that will be our main tools
in this paper.  Let $G_X(s;z,z')$ denote the Green's function (the integral kernel
of the resolvent operator $R_X(s) = (\Delta_X-s(1-s))^{-1}$).  
One way to cancel the kernel's singularity on
the diagonal so as to produce a trace is to subtract off the Green's function for
the model hyperbolic space, which has the same diagonal singularity.  This 
technique was used prominently by Patterson \cite{Patterson:1989}.
By lifting $G_X$ to $\bbH \times\bbH$, we can define
$$
\varphi_X(s;z) = (2s-1) \Bigl[ G_X(s;z,w) - G_{\bbH}(s;z,w) \Bigr]_{w=z},
$$
as a function on $\bbH$ (meromorphic in $s$).  By the invariance 
properties of the model hyperbolic Green's function it is easy to see that
$\varphi_X(s)$ descends to define a smooth function on $X$.   
This function is known to be polyhomogeneous in $\rho$
by the general theory of
\cite{MM:1987, Mazzeo:1991a}.  Or, more directly, one can deduce this from 
Guillop\'e's parametrix construction in \S1 of \cite{Guillope:1992}. 
This construction shows in particular that in the funnel ends, 
\begin{equation}  \label{infty.asymp1}
    \varphi_X(s)|_{F_j} \in \rho^{2s} \cinf(\bar F_j),
\end{equation} 
where $\bar F_j$ is the compactification of $F_j$ to the cylinder
$[0,1]_\rho\times (\bbR\bbZ)_x$.  
In the cusps,
\begin{equation}  \label{infty.asymp2}
   \varphi_X(s)|_{C_j} \in \rho^{2s-2} \cinf([0,1]_\rho) 
+ \rho^{-1} \cinf([0,1]_\rho) + O(\rho^\infty).
\end{equation}
(Coefficients of the asymptotic expansions in the cusps do not depend 
on the $x$ variable.)

It follows from (\ref{infty.asymp1}) and (\ref{infty.asymp2}) that, 
as $\varepsilon\to0$, 
$$
\int_{\{\rho>\varepsilon\}} \varphi_X(s;z)\>dg(z) \sim
a(s) + b_0(s) \log\varepsilon + \sum_{k=1}^\infty  b_k(s) \varepsilon^{k} +
\sum_{l=0}^\infty c_l(s) \varepsilon^{2s-1+l},
$$
and all of the coefficients are meromorphic by the meromorphy of the resolvents.
For  $s \ne 1/2-\bbN_0/2$ (where $\bbN_0 = \{0\}\cup\bbN$) we define
the function 
\begin{equation}\label{Phi.def}
\Phi_X(s) = \Oint_X \varphi_X(s;z)\>dg(z) = a(s),
\end{equation}
which extends meromorphically to the whole plane because $a(s)$ is meromorphic.
(The formula (\ref{Phi.def}) does not apply at the points $s\in 1/2-\bbN_0/2$ 
because the coefficient $c_l$ would also contribute to the 0-integral for 
$s = (1-l)/2$.)

Another function of interest is the renormalized trace of the spectral density:
\begin{equation}\label{Ups.def}
    \Upsilon_X(s) = (2s-1) \otr \Bigl[ R_X(s) - R_X(1-s) \Bigr].
\end{equation}
This 0-trace is well defined because the singularities of $G_X(s,z,w)$ and 
$G_X(1-s,z,w)$ on the diagonal cancel (and again because these kernels are 
polyhomogeneous in $\rho$).   
By examining the boundary asymptotics as above, we can see that $\Upsilon_X(s)$ 
also extends meromorphically to the whole plane, 
although the function is only given by the formula (\ref{Ups.def}) for $s \notin 
\frac12\bbZ$.

As one would expect, $\Phi_X(s)$ and $\Upsilon_X(s)$ are closely related. 

\begin{proposition} \label{UpPhiRelation}
 With $\Phi$ and $\Upsilon$ as defined above,
\begin{equation}\label{phi.ups}
        \Phi_X(s)~ +~ \Phi_X(1-s)~ =~ 
         \Upsilon_X(s)~ +~ (2s-1)\pi \cdot \chi(X) \cot(\pi s).
\end{equation}
\end{proposition}
\begin{proof}
For $s\notin \frac12\bbZ$, we have
\begin{equation*}
\begin{split}
\Phi_X(s) + \Phi_X(1-s) &= (2s-1) \Oint_X \Bigl[ G_X(s;z,w) - G_{X}(1-s;z,w) \\
&\qquad- G_{\bbH}(s;z,w) + G_{\bbH}(1-s;z,w)\Bigr]_{w=z} \>dg(z).
\end{split}
\end{equation*}
The difference $G_{\bbH}(s;z,w)-G_{\bbH}(1-s;z,w)$ is
continuous, and its restriction to the diagonal is necessarily a constant 
since these kernels depend only on $d(z,w)$.  This constant is shown
to be $(1/2) \cot \pi s$ in (\ref{ghgh.diag}), and
the result then follows from Lemma \ref{ovol.lemma}.
\end{proof}

The Selberg zeta function for the cylinder $M = \Gamma_\ell \backslash\bbH$ 
is
\begin{equation*}
              Z_{M}(s)~  =~ \prod_{k\geq 0}(1-e^{-(s+k)\ell })^2.
\end{equation*}
In Proposition 3.3 of \cite{Patterson:1989} Patterson computed 
\begin{equation}\label{patt.form}
\Phi_M(s) = \frac{Z_{M}'}{Z_{M}}(s).
\end{equation}
where $\Phi_M(s)$ is the function defined by (\ref{Phi.def}).

Now let $\Delta _{F}$ be the Laplacian on $F$ with Dirichlet boundary
conditions on the single closed geodesic bounding $F$.  We'll define 
\begin{equation}  \label{Z.F}
        Z_{F}(s)=e^{-s\ell /4}\prod_{k\geq 0}(1-e^{-(s+2k+1)\ell })^{2},
\end{equation}
where $\ell$ is the length of the geodesic boundary of $F$. 

Technically, $\Phi_X(s)$ was defined by (\ref{Phi.def}) only for 
$X$ a complete manifold.  To extend this definition to the manifold with 
boundary, $F$, we simply apply the same definition using $G_F(s)$ 
the Green's function for $\Delta_F$ and taking a 0-integral over $F$ 
(with regularization
needed only at the open end).  Together with the definition (\ref{Z.F}),
this allows us to make a statement analogous to (\ref{patt.form}):
\begin{proposition}\label{funnel.phi}
For $\re(s)>1/2$, 
\begin{equation}
\Phi_F(s) = \frac{Z_{F}'}{Z_{F}}(s).
\label{eq.hcyl.trace}
\end{equation}
\end{proposition}
\noindent
The proof will be deferred to Appendix \ref{sp.fun.sec}.

Since $\chi(F) =0$, by Proposition \ref{UpPhiRelation} we have 
\begin{equation}\label{F.tr}
\Upsilon_F(s) = \frac{d}{ds} \left( \frac{Z_{F}(s)}{Z_{F}(1-s)} \right).
\end{equation}

By Fourier analysis in the $x$ variable (in the coordinates 
(\ref{eq.Y.metric})),
one can write $\Delta _{F}$ as an infinite direct sum of ordinary
differential operators on the half-line; see \cite{GZ:1995a} for details.
This leads to the following explicit computation of the resonances of 
$\Delta _{F}$ as the set
\begin{equation} \label{PolarSet}
\left\{ \zeta _{k,m}=-(2k+1)+\frac{2\pi im}{\ell } \text{ (with multiplicity 2)}
:k\in\mathbb{N}_0,\;m\in \mathbb{Z}\right\} ,
\end{equation}
where $\ell = \ell(\partial F)$.
We can form a convergent Hadamard product  whose zero
set is (\ref{PolarSet}):
\begin{equation}
P_{F}(s)~  =~ \prod_{k=0}^{\infty }\prod_{m=-\infty }^{\infty }
(1-s/\zeta_{k,m})^2 e^{(2s/\zeta _{k,m}+(s/\zeta _{k,m})^{2})},
 \label{eq.PY}
\end{equation}
defining an entire function of order two. 
Note that (\ref{PolarSet}) is also the zero set of $Z_F$ defined in (\ref{Z.F}),
and hence
\begin{equation} \label{eq.YPZ}
        Z_{F}(s)~ =~ e^{p(s)} \cdot P_{F}(s)  
\end{equation}
for a polynomial $p$ of order at 
most two whose coefficients depend only on  $\ell$.

The calculation of $\Phi_H(s)$ for the horn $H$ is another
special function calculation which we will defer to Appendix 
\ref{sp.fun.sec}.  The result is:
\begin{proposition}\label{horn.phi}
For $H =  \Gamma_\infty\backslash  \bbH$ we have
\begin{equation}
\Phi_H(s) = -\log 2 - \Psi(s+1/2) + \frac1{2s-1},
\label{eq.horn.trace}
\end{equation}
where $\Psi(z)$ is the digamma function $\Gamma'(z)/\Gamma(z)$.
\end{proposition}
Because $\chi(H) = 0$, Propositions \ref{UpPhiRelation} and \ref{horn.phi}
immediately yield:
\begin{corollary}\label{horn.trace}
For $H =  \Gamma_\infty\backslash  \bbH$ we have
\begin{equation*}
\Upsilon_H(s) = -\frac{d}{ds}\log \left[4^s  
\frac{\Gamma(s-1/2)}{\Gamma(1/2 -s)} \right].
\end{equation*}
\end{corollary}


\section{Scattering theory}\label{scatt.sec}

As above, $X$ denotes a geometrically finite hyperbolic surface.  
We assume that $X$ has infinite area ($n_F\ge 1$), which implies that 
the positive Laplacian, denoted $\Delta _{X}$, has at most finitely
many eigenvalues in the interval $[0,1/4]$ and purely continuous spectrum in 
$[1/4,\infty )$ (no embedded eigenvalues). The resolvent 
$R_X(s)=(\Delta _{X}-s(1-s))^{-1}$, initially defined as a meromorphic 
$L^{2}(X) $-operator-valued function in the half-plane $\re (s)>1/2$, is
known to extend, as an operator from $C_{0}^{\infty }(X)$
to $C^{\infty }(X)$, to a meromorphic function in 
$\mathbb{C}$ (see Mazzeo-Melrose \cite{MM:1987} if $n_{C}=0$ and 
\cite{GZ:1997} if $n_{C}\neq 0$).   A pole of $R_X(s)$ is called a \textit{resonance}.
We let $\mathcal{R}_{X}$ denote the set of resonances, counted with
multiplicity $m_{\zeta }$ given for $\zeta\ne 1/2$ by the rank of the residue of
$R_X(s)$:
\begin{equation*}
m_\zeta = \rank \int_{\gamma _{\zeta }}(2s-1)R_{X}(s)\,ds,
\end{equation*}
where $\gamma _{\zeta }$ is a simple closed contour surrounding $\zeta$ and
no other resonance of $\Delta _{X}$.  (When $\zeta(1-\zeta)$ is an eigenvalue,
$m_\zeta$ is equal to the multiplicity.)  

The point $\zeta=1/2$ is exceptional because $R_X(s)$ has possibly a
second order singularity there.  For $\gamma_{1/2}$ a contour surrounding
$1/2$ and no other resonance, let
$$
a_{1/2} = \rank \int_{\gamma _{1/2}}(2s-1)R_{X}(s)\,ds,
$$
and
$$
b_{1/2} = \rank \int_{\gamma _{1/2}} R_{X}(s)\,ds.
$$
Lemma 4.1 of \cite{GZ:1997} shows that $a_{1/2}$ is the multiplicity of $1/4$ 
as an $L^2$-eigenvalue of $\Delta_X$ (which is zero if $n_C=0$).
To get the proper zero of the zeta function we set
$$
m_{1/2} = 2a_{1/2} + b_{1/2}.
$$
See the proof of Proposition \ref{Z.hlfplane} for the justification of this.

Let us recall some of the basics of scattering theory in our context; 
see \cite{GZ:1997} for details and complete proofs. Note that in \cite
{GZ:1997} the basic objects of scattering theory---the resolvent, generalized
eigenfunctions, and scattering operator---are defined as sections of certain
bundles over $X$; we in effect trivialize these bundles by our explicit
choice of defining function for the boundary of $X$ (see the remarks at the
beginning of \S \ref{sec.model}).

First we define the scattering operator $S_{X}(s)$ as a map from 
$\cinf(\bX) =  \mathbb{C}^{n_{C}}\oplus (\oplus
_{j=1}^{n_{F}}C^{\infty }(S^{1}))$ to itself.  
The kernel of $S_X(s)$ could be defined purely algebraically (as an average over
$\Gamma$ where $X = \Gamma\backslash\bbH$).  But a more intuitive
definition, following the general philosophy explained for example in
\cite{Melrose:1995}, is based on generalized eigenfunctions:

\begin{proposition}
Given $f = \{ z_{1},\cdots ,z_{n_{C}},f_{1},\cdots ,f_{n_{F}}\}\in C^{\infty }(\bX)$, 
$s\notin \bbZ/2$, there is a unique
smooth solution $u$ of the eigenvalue problem 
\begin{equation*}
\left( \Delta _{X}-s(1-s)\right) u=0
\end{equation*}
with 
$$
u\sim \rho ^{s}\tilde f+\rho^{1-s}\tilde g,
$$
where $\tilde f, \tilde g \in\cinf(X)$ are polyhomogeneous of degree zero in $\rho$,
and $\lim_{\rho\to 0} \tilde f = f$.
\end{proposition}

See \cite{JS:1998} for a proof which covers our case $n_C=0$;
the generalization to $n_C>0$ is straightforward using the parametrix construction
from \cite{Guillope:1992}.
From this proposition we obtain a linear map $S_X(s): 
f\mapsto g =  \lim_{\rho\to0} \tilde g$. 
It follows from the definition that 
$S_{X}(s)S_{X}(1-s)$ is the identity operator.

In what follows, we write $i$ or $i'$ for a cusp index $1\leq
i,i'\leq n_{C}$ and $j$ or $j'$ for a funnel index $
1\leq j,j'\leq n_{F}$. We denote by $S_X^{ii'}(s)$ the
component of $S_{X}(s)$ mapping $\partial_\infty C_{i'}$ to 
$\partial_\infty C_{i}$, etc.  From the detailed description of
$S_X(s)$ in \cite{GZ:1997} we extract:

\begin{proposition}
The map $S_{X}(s)$ extends to a meromorphic operator-valued
function with the following properties:
\begin{enumerate}
\item  $\left\{ S_X^{ii'}(s):1\leq i,i'\leq n_{C}\right\} $ is
a matrix-valued meromorphic function of $s$,
\item  $S_X^{jj'}(s)$ is a smoothing operator if $j\neq j'$,
and $S_X^{jj}(s)-S_{F_{j}}(s)$ extends to a meromorphic family of smoothing
operators on $C^{\infty }(S^{1})$, and 
\item 
$S_X^{ij}(s)$ and $S_X^{ji}(s)$ have integral kernels in $C
^{\infty }(S^{1})$ which are meromorphic functions of $s$.
\end{enumerate}
\end{proposition}

The poles of the scattering operator can be deduced from the resonance set 
$\mathcal{R}_X$, as shown in \S2 of \cite{GZ:1997}.
(They coincide precisely except for infinite rank poles of $S_X(s)$ at positive 
integer points, and possibly finitely many points related to the discrete spectrum.)  

In \cite{GZ:1995a}, Guillop\'{e} and Zworski proved that number of resonances 
(including multiplicities) in the ball $|\zeta|<r$ 
grows at most quadratically in $r$:
\begin{equation*}
     \sum_{\zeta \in   \mathcal{R}, |\zeta|<r}  m_{\zeta} ~ =~ O( r^2).
\end{equation*}
Hence the Hadamard product 
\begin{equation}\label{px.def}
     P_{X}(s)~ =~ s^{m_{0}}\prod_{\zeta \in \mathcal{R}}
(1-s/\zeta)^{m_{\zeta }}e^{m_{\zeta }(s/\zeta +s^{2}/2\zeta ^{2})}
\end{equation}
converges.

We recall from \cite{GZ:1997} that the relative scattering determinant 
$\tau _{X}(s)=\det (S_{rel}(s))$, defined in \S\ref{intro.sec}, 
is determined up to finitely many parameters by the
lengths $\ell _{j}$ of geodesics bounding the $F_{j}$ together with the
resonance set $\mathcal{R}_{X}$. To state this result more precisely, recall
the definition (\ref{eq.PY}) of $P_{F}(s)$ and define
\begin{equation*}
P_{Y}(s)=\prod_{j=1}^{n_{F}}   P_{F_{j}}(s).
\end{equation*}
\begin{proposition}[Proposition 3.4 of \cite{GZ:1997}] \label{prop.tau1}
The function $\tau _{X}$ extends to a meromorphic function
in the complex plane having the form
\begin{equation*}
     \tau _{X}(s)~ =~ e^{h(s)}\frac{P_{X}(1-s)}{P_{X}(s)}
      \frac{P_{Y}(s)}{P_{Y}(1-s)}
\end{equation*}
where $h(s)$ is a polynomial of degree at most four.
\end{proposition}
\noindent
(We will find that $h(s)$ has degree at most two at the end of 
\S\ref{zeta.sec}.)

Let $R_{X}(s)=\left( \Delta _{X}-s(1-s\right) )^{-1}$ and let 
$R_{Y}(s)=(\Delta _{Y}-s(1-s))^{-1}$, where $Y$ is the union of the funnels
as before with Dirichlet conditions on $\partial Y$.  Note that the operators 
$R_{X}(s)-R_{X}(1-s)$ and $R_{Y}(s)-R_{Y}(1-s)$ actually have smooth
kernels although with insufficient decay at infinity to be trace-class.
There is a natural identification of $Y$ with the corresponding
submanifold of $X$ and we denote by $\mathbf{1}_{Y}$ the obvious induced
mapping from $C^{\infty }(X)$ onto $C^{\infty }(Y)$.
We set
\begin{equation*}
    Q_{X}(s)=R_{X}(s)-\mathbf{1}_{Y}^{t} \circ R_{Y}(s) \circ \mathbf{1}_{Y}\text{.}
\end{equation*}
\begin{proposition}[Equation (4.14) of \cite{GZ:1997}\footnote{We believe there
to be a typo in formula (4.14) of \cite{GZ:1997}: the terms on the 
right-hand side of their equation should carry negative signs.}]
\label{trace.prop}
For $\re(s) = \frac12$, $s\ne 1/2$,
\begin{equation} \label{eq.tr1}
(2s-1)\cdot  \otr \left[  Q_{X}(s)-Q_{X}(1-s)\right]~ 
=~ - \frac{\tau _{X}'}{\tau _{X}}(s).
\end{equation}
\end{proposition}
In the original formula in \cite{GZ:1997}, the trace in (\ref{eq.tr1}) is interpreted
as a distribution on the line $\re(s) = 1/2$ and the right-hand side includes 
a delta-function singularity at $s=1/2$.  The singularity occurs for reasons discussed
in \S\ref{sec.model}; the same issue is behind our restriction of 
formula (\ref{Ups.def}) to $s\notin \frac12\bbZ$.


\section{Selberg's zeta function}\label{zeta.sec}

Let $X = \Gamma\backslash\bbH$ be a geometrically finite hyperbolic surface.
Recall that Selberg's zeta function, $Z_X$, is defined 
for $\re (s)>1$ by
\begin{equation}  \label{eq.zeta.def}
     Z_{X}(s)=\prod_{\left\{ \gamma \right\} }\prod_{k=0}^{\infty }
     \left( 1-e^{-(s+k)\ell (\gamma )}\right). 
\end{equation}
where the outer product goes over conjugacy classes of primitive hyperbolic
elements of $\Gamma$, and $\ell (\gamma )$ is the length of the corresponding 
closed geodesic.  Following \cite{Sarnak:1987}, we define  
\begin{equation*}
      Z_{\infty }(s)~ =~ 
       \left[\frac{(2\pi )^{s}\Gamma _{2}(s)^{2}}{\Gamma (s)}\right]^{-\chi}.
\end{equation*}
where $\Gamma _{2}(s)$ is Barnes' double Gamma function
and $\chi= \chi(X)$ is the Euler characteristic of $X$. 
The function $Z_\infty(s)$ has no zeros, and poles of order
$(2k+1)(-\chi)$ at  $s= -k$, $k\in\bbN_0$.

The purpose to this section is to prove the following
result. 
\begin{theorem} \label{thm.zeta2}
The Selberg zeta function has a factorization
\begin{equation}   \label{zetafact}
      Z_{X}(s) \cdot Z_{\infty }(s)~  =~ e^{q(s)} \cdot 
     \Gamma (s-1/2)^{n_{C}} \cdot P_{X}(s),
\end{equation}
where $q(s)$ is a polynomial of degree $\leq 2$ and $P_X(s)$ is the Hadamard
product over resonances (\ref{px.def}).
\end{theorem}

Theorem \ref{thm.zeta} is an immediate consequence of Theorem \ref{thm.zeta2}.
But before taking up the proof of Theorem \ref{thm.zeta2}, we will show how
Corollary \ref{cor.finite.euler} follows.
\begin{proof}[Proof of  Corollary \ref{cor.finite.euler}]
Suppose that the resonance set is fixed, which determines $P_X(s)$.
We claim that $P_X(s)$ determines $\chi$ (hence $Z_\infty$), $n_C$, and $q(s)$,
and therefore fixes $Z_X(s)$, from which the length spectrum may be
deduced by a standard argument.
To see this we take the log of (\ref{zetafact}) and analyze the asymptotics as 
$\re (s)\rightarrow \infty$.  On the left-hand side, $\log Z_X(s)$ decays 
exponentially, while 
$$
\frac1{\chi}\log Z_\infty(s) \sim \frac{2+\log2\pi}{2}
-2\zeta^{\prime}(-1)-\Bigl(\frac{1}{2}s(s-1)-\frac{1}{6}\Bigr)\log s(s-1)
+\frac{3}{2}s(s-1).
$$
By Stirling's formula,
$$
\log \Gamma(s-1/2) \sim \frac{1+\log2\pi}{2} + s\log (s-1/2) - s.
$$
This implies that $\log P_X(s)$ has an asymptotic expansion as $\re(s)\to\infty$. 
The value of $\chi$ may be read off from the $s^2 \log s$ coefficient. 
Then, after subtracting off the $\log Z_\infty(s)$ term, the $s\log s$ coefficient
determines $n_C$.   Once the $n_C$ terms are removed, what remains in
the asymptotic expansion is precisely $q(s)$, so this is also determined by $P_X(s)$.

Now assume that the length spectrum is known, giving $Z_X(s)$.  
The order of the zero of $Z_X(s)$ at $s=-k$, $k\in\bbN_0$, is $m_{-k} 
+ (2k+1)(-\chi)$.  This gives us for any $k$ a bound 
$$
0\ge \chi \ge \frac{1}{2k+1}\ord_{s=-k} Z_X(s),
$$ 
implying that only finitely many values of $\chi$ are possible.  And hence finitely many 
choices of $n_C$, since $n_C \le 2 - \chi$.  If $n_C$ and $\chi$ were fixed,
then $Z_X(s)$ would determine the divisor of $P_X(s)$ and the resonance 
set could be read off  directly.
\end{proof}

\begin{remark}
Corollary \ref{cor.finite.euler} could also have been deduced from the 
wave trace formula
proved by Guillop\'e and Zworski in \cite{GZ:1999}.  The combination of
Theorem 1 and Equation (3.1) of that paper yields the following trace formula. 
As distributions on $\bbR_+$,
\begin{equation*}\begin{split}
\sum_{\zeta\in\mathcal{R}_X} m_\zeta e^{(\zeta-1/2)t} &=  \sum_{\{\gamma\}}
\sum_{k=1} \frac{\ell(\gamma)}{\sinh(k\ell(\gamma)/2)} \delta(t - k\ell(\gamma)) \\
&\qquad +\chi \frac{\cosh(t/2)}{2\sinh^2(t/2)}
+ \frac{n_C}{2} \coth(t/2),
\end{split}\end{equation*}
where the $\{\gamma\}$ denotes a list of conjugacy classes of primitive hyperbolic
elements of $\Gamma$, as in the definition of the zeta function.
Note that 
$$
\frac{\cosh(t/2)}{2\sinh^2(t/2)} = \sum_{k=0}^\infty (2k+1) e^{(-k-1/2)t},
$$
for $t>0$.   The problem of overlap between resonances and topological zeros would
not be avoided by this alternate route.
\end{remark}

\begin{proposition}\label{patt.formula}
The logarithmic derivative of the zeta function is given by
\begin{equation}\label{pform.eq}
     \frac{Z_X'}{Z_X}(s) = \Phi_X(s) - n_C \Phi_H(s) 
\end{equation}
\end{proposition}
\begin{proof}
Since both sides of (\ref{pform.eq}) are known to be meromorphic, it suffices to 
consider $\re(s)>1$.

For $\re(s)>1$, the Green's function $G_X(s;z,w)$ can be written as a 
convergent sum
\begin{equation*}
\begin{split}
      G_X(s;z,w) &= G_\bbH(s;z,w) + \sum_{\gamma\in \Pi_h} 
      \sum_{g \in \Gamma/\langle  \gamma\rangle}
      \sum_{k\ne 0} G_\bbH(s;z,g^{-1}\gamma g w) \\
  &\qquad +  \sum_{\gamma\in \Pi_p} \sum_{g \in \Gamma/\langle \gamma\rangle}
\sum_{k\ne 0} G_\bbH(s;z,g^{-1}\gamma g w),
\end{split}\end{equation*}
where $\Pi_h$ and $\Pi_p$ are lists of representatives of conjugacy 
classes of maximal hyperbolic and parabolic subgroups of $\Gamma$, 
respectively.  It follows that 
$$
\Phi_X(s) = (2s-1) \int_X [G_X(s;z,w) - G_\bbH(s;z,w)]_{z=w} \:dg(z) 
$$
into corresponding sums over $\Pi_h$ and $\Pi_p$.  
For any $\gamma\in \Gamma$, the standard
trace formula technique for summing over a conjugacy class gives:
\begin{equation*}\begin{split}
\sum_{g \in \Gamma/\langle \gamma\rangle} \sum_{k\ne 0}  \int_X
G_\bbH(s;z,g^{-1}\gamma g w) &= \int_{\langle\gamma\rangle\backslash\bbH}
[G_{\langle\gamma\rangle\backslash\bbH}(s;z,w) - G_\bbH(s;z,w)]_{z=w} \:dg(z)\\
&= \Phi_{\langle\gamma\rangle\backslash \bbH}.
\end{split}\end{equation*}
If $\gamma$ is hyperbolic, then by Patterson's formula 
(\ref{patt.form}) this is equal to 
$(\log Z_{\langle\gamma\rangle\backslash \bbH})'(s)$.  
And then since, by definition, 
$$
Z_X(s) = \prod_{\gamma\in\Pi_h} Z_{\langle\gamma\rangle\backslash \bbH}(s),
$$
summing these hyperbolic terms over $\Pi_h$ 
gives the logarithmic derivative of $Z_X(s)$.

For any $\gamma$ parabolic, the quotient $\langle\gamma\rangle\backslash\bbH$ 
is isometric to the model horn $H$ from \S\ref{sec.model}, so
the contribution from each cusp is $\Phi_H(s)$.
\end{proof}

Since $Y$ is a direct sum of the funnels $F_j$, we set
$$
Z_{Y}(s)=\prod_{j=1}^{m}Z_{F_{j}}(s) = \prod_{j=1}^{m}
\left[ e^{-s\ell_j /4}\prod_{k\geq 0}(1-e^{-(s+2k+1)\ell_j })^{2}\right],
$$

\begin{proposition}\footnote{In the case $n_C = 0$, such a formula was noted in \cite{GZ:1997}, equation (5.3).}
\label{prop.tau2}
For some constant $c$,
\begin{equation*}
\tau _{X}(s) = e^{c - \left(n_{C}\log 4\right) s}\frac{Z_{X}(1-s)}{Z_{X}(s)}
\frac{Z_{\infty }(1-s)}{Z_{\infty}(s)}
\frac{Z_{Y}(s)}{Z_{Y}(1-s)}\left( \frac{\Gamma (s-1/2)}{\Gamma (1/2-s)}
\right) ^{n_{C}}.
\end{equation*}
\end{proposition}
\begin{proof}
From Proposition \ref{trace.prop} and (\ref{Ups.def}) we find that
for  $s\notin \frac12\bbZ$ 
$$
    -\frac{\tau_X'(s)}{\tau_X(s)}~  =~  \Upsilon_X(s) - (2s-1) \otr_X 
\Bigl[\mathbf{1}_{Y}^{t} \circ \bigl(R_{Y}(s) - R_{Y}(1-s)\bigr) 
\circ \mathbf{1}_{Y}\Bigr].
$$
Since $R_{Y}$ breaks up into a direct sum of $R_{F_j}$ on the
individual funnels, the second term on the right-hand side can be reduced to
$$
(2s-1) \otr_X \Bigl[\mathbf{1}_{Y}^{t} \circ \bigl(R_{Y}(s) - R_{Y}(1-s)\bigr) \Bigr]
= \sum_{j=1}^{n_F} \Upsilon_{F_j}(s),
$$
by the definition of $\Upsilon_{F_j}(s)$.
From (\ref{F.tr}) we find
\begin{equation}
\sum_{j=1}^{n_F} \Upsilon_{F_j}(s) = \frac{d}{ds}
\log \left( \frac{Z_{Y}(s)}{Z_{Y}(1-s)}\right).  \label{Y.tr}
\end{equation}

Given the computation of $\Upsilon_H(s)$ in Corollary \ref{horn.trace}, 
to prove the Proposition it will suffice to show that
\begin{equation}  \label{upsfunc.eqn}
   \Upsilon_X(s) =\frac{d}{ds}\log 
    \left( \frac{Z_{X}(s)}{Z_{X}(1-s)}\frac{Z_{\infty}(s)}{Z_{\infty}(1-s)}\right) 
   +n_{C}\Upsilon_H(s).
\end{equation}
To deduce (\ref{upsfunc.eqn}), one can use
the product expansion of $\Gamma_2(s)$ to 
derive the identity
\begin{equation}     \label{zinf.id}
    \frac{Z_{\infty }'(s)}{Z_{\infty }(s)}~  =~ \chi \cdot (2s-1) \cdot (\Psi (s)-1).
\end{equation}
where $\Psi (s)=\Gamma'(s)/\Gamma (s)$.  
It then follows from the functional equation 
$\Gamma(s) \cdot \Gamma (1-s)=\pi \csc (\pi s)$ that
\begin{equation}\label{zinf.eq}
      \frac{d}{ds}\log \left( \frac{Z_{\infty }(s)}{Z_{\infty }(1-s)}\right)~
      =~ - \chi \cdot (2s-1) \cdot \pi \cdot \cot (\pi s)
\end{equation}
Equation (\ref{upsfunc.eqn}) follows by combining  equation (\ref{zinf.eq}),  
Proposition \ref{UpPhiRelation},  and     Proposition \ref{patt.formula}.
\end{proof}

\begin{proposition}\label{Z.hlfplane}
The meromorphic function $Z'_X(s)/Z_X(s)$ is analytic in $\re(s)>1/2$ except for 
first-order poles at points $\zeta$ where $\zeta(1-\zeta)$ is an 
eigenvalue (with residues equal to the eigenvalue multiplicities).  Also, 
$Z'_X(s)/Z_X(s)$ has no poles on $\re(s) = 1/2$ except possibly at $s=1/2$.
The residue of $Z'_X(s)/Z_X(s)$ at $s=1/2$ is $m_{1/2} - n_C$.
\end{proposition}
\begin{proof}
Except at $s=1/2$ the argument is essentially the same as in Theorem 6.2 of 
\cite{PP:2000}.  Recall the integrand
$\varphi_X(s; w)$  appearing in the definition (\ref{Phi.def}) of $\Phi_X(s)$. 
The poles of
$\varphi_X(s)$ for $\re(s)\ge 1/2$, $s\ne 1/2$, and its behavior near those poles 
can be deduced from Lemmas 2.4 of \cite{GZ:1997}.  In particular, in $\re(s)>1/2$,
$\varphi_X(s)$ is meromorphic with poles only at points $\zeta \in (1/2,1)$
for which $\zeta(1-\zeta) \in \sigma(\Delta_X).$  Near such a point $\zeta$ 
it has the form
$$
\varphi_X(s) = \frac1{s-\zeta} \sum_{k=1}^{m_\zeta} \psi_k^2 + \tilde\varphi(s),
$$
where $\{\psi_k\}$ form an orthonormal basis of the eigenspace for $\zeta(1-\zeta)$,
and $\tilde\varphi(s)$ is holomorphic near $s=\zeta$.  Since $\tilde\varphi(s)$ will
be integrable for $\re(s)>1/2$, from Proposition \ref{patt.formula}
we have  near $s=\zeta$:
$$
\frac{Z_X'}{Z_X}(s) = \frac{m_\zeta}{s-\zeta} + h(s).
$$
where $h$ is holomorphic.

For $\re(s) = 1/2$, $s \ne 1/2$, the function $\varphi_X(s)$ is holomorphic
(there are no embedded eigenvalues as $n_F>0$),  
but it is not integrable.  Nonetheless, from the behavior of $\varphi_X(s)$
at infinity  (see (\ref{infty.asymp1}) and  (\ref{infty.asymp2})) it
follows that the $0$-integral of $\varphi_X(s)$ vanishes for this case.

This leaves finally the point $s = 1/2$.  
By Lemma 4.1 of \cite{GZ:1997} the structure of $R_X(s)$ near $s=1/2$ is
\begin{equation}\label{rabc}
R_X(s) = \frac{A}{(2s-1)^2} + \frac{B}{(2s-1)} + C(s),
\end{equation}
where $C(s)$ is analytic near $s=1/2$.  Our definitions in \S\ref{scatt.sec}
set $a_{1/2} = \rank A$,  $b_{1/2} = \rank B$, and $m_{1/2} = 2a_{1/2} + b_{1/2}$.
Substituting (\ref{rabc}) into the definition of $\varphi_X$ shows that
$$
\varphi_X(s;z) = \frac{a(z)}{2s-1} + b(s;z),
$$
where $b(s;z)$ is analytic near $s=1/2$.

Suppose $u$ is an eigenvalue of $\Delta_X$ with eigenvalue $1/4$. 
Applying $R_X(s)$ to
$$
(\Delta_X - s(1-s))u = (s-1/2)^2 u
$$
shows that $u = \frac14 Au$.  Hence $A/4$ is the projector onto the 
$1/4$-eigenspace of $\Delta_X$.
The contribution of the $a$ term to the residue of $\Phi_X(s)$ at $s=1/2$ is 
$$
\frac12 \int_X a \>dg = \frac12\tr A = 2 a_{1/2}.
$$

Although $b(s;z)$ is analytic near $s=1/2$, its 0-integral does contribute
to the pole in $\Phi_X(s)$ because of its boundary asymptotics.   In a funnel $F_j$,
$b(s;\cdot)|_{F_j} \in \rho^{2s}\cinf{\bar F_j}$ and the contribution can be
deduced exactly as in  Theorem 6.2 of \cite{PP:2000}.   Namely, 
the 0-integral of $b|_{F_j}$ 
contributes a residue of $\frac12 \tr[S_X^{jj}(1/2) + I]$.  (For comparison of 
formulas, note that the normalized scattering matrix $\mathcal{S}_X(s)$ used in 
\cite{PP:2000} is related to our definition by $\mathcal{S}_X(1/2)= - S_X(1/2)$.)

The contribution of a cusp is slightly different.  Consider cusp $C_i$ and introduce 
a cutoff $\eta$ which is supported in $C_i$ and equal to 1 in a neighborhood of 
$\rho=0$.  The contribution from cusp $C_j$ is
$$
\res_{s=1/2} \Oint_{C_j} \eta(z) \phi_X(s;z) \>dg(z)
$$
(which will independent of $\chi$.)  Identifying $C_j$ with the cusp end of $H$,
and using the parametrix construction in \cite{Guillope:1992} we can identify 
$$
\varphi_X(s;z)|_{C_i} = \varphi_H(s;z) + \rho^{2s-2} k(s;z),
$$
where, under our convention for $S_X(s)$, $\lim_{\rho\to0} k(s,z)  
= S_X^{ii}(s)$.  The contribution
to the residue from the $k(s;z)$ term is then easily computed to be 
$\frac12 S_X^{ii}(1/2)$. Lemma \ref{cusp.res} shows that
$$
\res_{s=1/2} \Oint_{C_j} \eta(z) \varphi_H(s;z) \>dg(z) = 0,
$$
so the $\varphi_H(s;z)$ term does not contribute to the residue.

The full contribution of the $b$ term to the residue of $\Phi_X(s)$
at $s=1/2$ is thus $(1/2) \tr[S_X^{FF}(s) + I] + (1/2) \tr[S_X^{CC}(s)]$.
By adding and subtracting $\frac12 \tr I^{CC} = n_C/2$ we can write
this as $\tr[S_X(1/2) + I] - n_C/2$.  By
Lemma 4.3 of \cite{GZ:1997} (with our convention for $S_X(s)$),
$$
\frac12 \tr[S_X(1/2) + I] = b_{1/2}.
$$
Thus the contribution of the $b(s;z)$ term to the residue is
$b_{1/2} - n_C/2$.

Finally, combining the $a$ and $b$ contributions shows that
$$
\res_{s=1/2} \Phi_X(s) = 2a_{1/2} + b_{1/2} - n_C/2.
$$
And the stated result follows because $\res_{s=1/2}\Phi_H(s) = 1/2$.
\end{proof}

\begin{proof}[Proof of Theorem \ref{thm.zeta2}]
Combining Propositions \ref{prop.tau1} and \ref{prop.tau2} and 
equation (\ref{eq.YPZ}) yields
\begin{equation}   \label{eq.zz}
\frac{Z_{X}(s)Z_{\infty }(s)}{Z_{X}(1-s)Z_{\infty }(1-s)}~  =~ e^{h_1(s)}
\frac{P_{X}(s)}{P_{X}(1-s)} \cdot \left( \frac{\Gamma
(s-1/2)}{\Gamma (1/2-s)}\right) ^{n_{C}} 
\end{equation}
for a polynomial $h_1$ of order at most four. 

Now let
$$
G(s) = \frac{Z_X(s) Z_\infty(s)} {P_X(s) \Gamma(s-1/2)^{n_C}}.
$$  
Proposition \ref{Z.hlfplane} says that $G(s)$ has no zeroes or poles in 
$\re(s)\ge 1/2$.  And (\ref{eq.zz})
says $G(s) = e^{h_1(s)} G(1-s)$, implying that $G$ has no zeroes or
poles in $\re(s) \le 1/2$ either.  Hence $G$ is the exponential of an
entire function $q(s)$.

It remains to show that $q(s)$ is a polynomial 
with degree equal to 2.  From the convergent
Euler product (\ref{eq.zeta.def}) and elementary estimates on the counting
function for lengths of closed geodesics, we can easily show that $Z_{X}(s)$
is bounded in the half-plane $\re (s)\ge 2$.  Since $P_X(s)$ and $1/Z_\infty(s)$
are entire of order 2 and $1/\Gamma(s-1/2)$ is order 1, we conclude that
$|q(s)| \le C_\kappa |s|^{2+\kappa}$ for $\kappa>0$ and $\re(s)\ge 2$. 
The functional
relation \ref{eq.zz} and the fact that $h_1$ is order 4 then imply a bound
$|q(s)| \le C|s|^{4}$ for $\re(s)\le -1$.

To handle the strip $S_\varepsilon = \{-1-\varepsilon< \re(s) < 2 + \varepsilon\}$, 
consider the function 
$$
F_X(s) =  \otr [Q_X(s) - Q_X(s_0)] 
$$
for $s_0$ fixed and not a pole of $Q_X(s)$, and similarly $F_H(s)$
for $H = \Gamma_\infty\backslash \bbH$.  By taking derivatives of 
the formula (\ref{pform.eq}), we can derive
$$
\left( \frac1{2s-1} \frac{d}{ds} \right)^2 \log \frac{Z_X(s)}{Z_Y(s)} = 
\left( \frac1{2s-1} \frac{d}{ds} \right) \Bigl[F_X(s) + \chi \Psi(s) - n_C  
F_H(s)\Bigr]
$$
Since $Z_Y$ and $\Psi$ are explicit, to control the growth of $\log Z_X(s)$
in $S_\varepsilon$ it suffices to control $F_X(s)$ (and $F_H(s)$).

Introduce, as in \S5 of \cite{BJP:2002}, disks $D_j$ such that $\mathcal{R}_X \subset 
\cup_j D_j$ and $d(s, \mathcal{R}_X) \ge C\langle s\rangle^{-2-\delta}$ for all
$s\in \bbC\setminus(\cup_j D_j)$.  Using the parametrix
construction of \cite{GZ:1997} and arguing exactly as in Appendix B of
\cite{BJP:2002}, we obtain the bound
$$
|F_X(s)| \le C(\eta) \exp(|s|^{2+\eta}),
$$
for $s \in S_\varepsilon\setminus(\cup_j D_j)$ (and we can similarly bound 
$F_H$).  Using the maximum modulus principle we obtain from these
estimates a bound $|q(s)| \le C(\eta) \exp(|s|^{2+\eta})$ for $s\in S_\varepsilon$.
The Phragm\'{e}n-Lindel\"{o}f theorem then allows us to extend the polynomial
bound on $q(s)$ into the strip $S_\varepsilon$ and
conclude that the $q(s)$ in (\ref{zetafact}) is a polynomial of degree at most four. 

To improve the order bound, we examine the log of (\ref{zetafact}) as 
$\re (s)\rightarrow \infty$.  In this limit, $\log Z_X(s)$ decays exponentially,
while $\log Z_\infty(s) = O(s^2\log s)$ and $n_C \log \Gamma(s-1/2) =
O(s\log s)$.   Since $P_X$ is entire of order 2, it follows that $q(s)$
has degree no more than two.  
\end{proof}

With Theorem \ref{thm.zeta2} proven, we can revisit Proposition
\ref{prop.tau2} and deduce the following:
\begin{corollary}
The relative scattering determinant $\tau_X(s)$ is the ratio of entire functions
of order two (i.e. the polynomial $h(s)$ in Proposition \ref{prop.tau1} has degree
at most two).
\end{corollary}


\section{Determinant of the Laplacian}

For $X$ an infinite area geometrically finite Riemann surface, we can 
define (as in \cite{BJP:2002}) a determinant of the Laplacian $D_X(s)$ formally equal
to $\det(\Delta - s(1-s))$ by integrating
\begin{equation}\label{det.def}
\left( \frac1{2s-1} \frac{d}{ds}\right)^2 \log D_X(s) = - \otr R_X(s)^2.
\end{equation}
(The right-hand side is well-defined because $R_X(s)^2$ has a continuous kernel.)
In this section we'll derive explicit formulas for $D_X(s)$ so defined.

Setting 
$$
\mathcal{L}_s = \frac{1}{2s-1} \frac{d}{ds},  
$$
for convenience, 
we begin by differentiating the definition of
$\Phi_X(s)$ and separating the traces:
$$
\mathcal{L}_s \left(\frac{\Phi_X(s)}{2s-1}\right) = - \otr R_X(s)^2 + \otr 
R_\bbH(s)^2.
$$
The second term on the right makes sense because the restriction
of the integral kernel of $R_{\bbH}(s)^2$ to the diagonal is a constant. In
(\ref{rsq.diag}) this constant is computed to be
$(2\pi)^{-1} \mathcal{L}_s \Psi(s)$.  Since $\ovol(X) = -2\pi\chi$, we obtain
$\otr  R_\bbH(s)^2 =  - \chi  \mathcal{L}_s \Psi(s)$, and by (\ref{zinf.id}) this 
is in turn equal to $-\mathcal{L}_s^2\log Z_\infty(s)$.  We thus have
$$
\mathcal{L}_s \left(\frac{\Phi_X(s)}{2s-1}\right) =
- \otr R_X(s)^2 - \mathcal{L}_s^2 \log Z_\infty(s).
$$

By differentiating (\ref{pform.eq}) we could also write
$$
\mathcal{L}_s \left(\frac{\Phi_X(s)}{2s-1}\right)  =
\mathcal{L}_s^2 \log Z_X(s)  - n_C
\mathcal{L}_s \left(\frac{\Phi_H(s)}{2s-1}\right),
$$
and by Proposition \ref{horn.phi}
$$
\frac{\Phi_H(s)}{2s-1} = - \mathcal{L}_s \log 
\left[2^s (s-1/2)^{1/2} \Gamma(s-1/2) \right].
$$
Putting these formulas together, we see that
\begin{equation}\label{rsq.trace}
- \otr R_X(s)^2 = \mathcal{L}_s^2 \log \left[ 
\frac{Z_X(s)Z_\infty(s)} {2^{s n_C} (s-1/2)^{n_C/2}\Gamma(s-1/2)^{n_C}}
\right]
\end{equation}
In view of Theorem \ref{thm.zeta2}, this proves:
\begin{theorem}
If $D_X(s)$ is defined by (\ref{det.def}), then $(s-1/2)^{n_C/2} D_X(s)$ is an entire
function with zeros given by $\mathcal{R}_X$ (with multiplicities),
satisfying
\begin{equation*}\begin{split}
(s-1/2)^{n_C/2} D_X(s) &= e^{[Fs(1-s) + E+ sn_c\log 2]} \frac{Z_X(s)Z_\infty(s)}
{\Gamma(s-1/2)^{n_C}} \\
&= e^{q_1(s)} P_X(s),
\end{split}\end{equation*}
where $F$, $E$ are constants and $q_1(s)$ is a polynomial of degree $\le 2$.
\end{theorem}

This is consistent with previous calculations of the determinant
in \cite{BJP:2002, Efrat:1988, Sarnak:1987}.
Another way to interpret (\ref{rsq.trace}) is as an analytic
expression for the second derivative of $\log P_X(s)$:
$$
\mathcal{L}_s^2 \log P_X(s) = - \otr R_X(s)^2 + \frac{A}{(2s-1)^3}
- \frac{4n_C}{(2s-1)^4},
$$
for some constant $A$.


\section{Isospectral finiteness for geometrically 
finite hyperbolic surfaces}

The purpose of this section is to prove Theorem \ref{Finiteness}.
First, we make the notion of {\em length spectrum counted with multiplicities}
more precise.  Given a complete hyperbolic surface $X$, let $N_X(L)$
denote the number of closed geodesics on $X$ with 
length less than $L$.   The function $N_X$ is the length spectrum 
counting function.  We say that two surfaces $X$ 
and $Y$ are {\em (length) isospectral (including multiplicities)} 
if and only if  $N_X= N_Y$ as functions. 

By Corollary \ref{cor.finite.euler} the resonance set determines the 
length spectrum, and the length spectrum fixes the Euler characteristic 
up to finitely many possibilities.
The Euler characteristic in turn determines the homeomorphism class 
up to finitely many possibilities.
So to prove Theorem \ref{Finiteness} it suffices to prove that if $X$ is a
geometrically finite hyperbolic surface of infinite area, then there 
are only finitely many hyperbolic surfaces $Y$, homeomorphic to $X$, 
such that $Y$ is (length) isospectral to $X$.  Moreover, as there are 
only finitely many ends, we may assume 
that the number of cusps and number of funnels are constant.

We restate and prove this claim
in the (hyperbolic) language of  Teichm\"uller theory
 \cite{Buser:1992, Trm}. Let $S$ be a fixed smooth 
surface  of finite Euler characteristic with $\partial S = \emptyset$. 
Here $g$ will denote a hyperbolic metric 
which is (Riemannian) complete and has exactly $n_C$ 
cusps and $n_F$ funnels. The set of all such metrics $g$ on 
$S$ has a natural topology coming from the $C^{k}$ convergence 
of tensors on $S$ \cite{Trm}. The {\em moduli space}, 
${\mathcal M}(S)$, is the quotient of this space by the group
of smooth self-diffeomorphisms of $S$ acting via pull-back. 
In other words,  ${\mathcal M}(S)$ is the set of isometry 
classes $[g]$ of complete hyperbolic metrics $g$ on $S$.

The {\em Teichm\"uller space} ${\mathcal T}(S)$ of $S$ is 
the quotient of the space of hyperbolic metrics by the 
diffeomorphisms isotopic to the identity. 
It is well-known that ${\mathcal T}(S)$ is naturally 
homeomorphic to a ball in some Euclidean space, and that the 
natural projection $\pi: {\mathcal T}(S) \rightarrow {\mathcal M}(S)$ 
is an orbifold covering map.

In this language, Theorem \ref{Finiteness} is reduced by virtue of
Corollary \ref{cor.finite.euler} to the following:

\begin{theorem}  \label{IsofiniteTeich}
Let ${\mathcal A} \subset {\mathcal M}(S)$ be an isospectral
set of isometry classes. Then ${\mathcal A}$ is a finite set.
\end{theorem}

To prove this theorem\footnote{During the preparation of this work, we received a
preprint from Inkang Kim \cite{Kim} with an isofiniteness result (for 
fixed homeomorphism class). We do not fully understand his argument.}, 
we first note that 
it suffices to prove that the isospectral set ${\mathcal A}$
is precompact.  Indeed, this would imply that the lift of 
${\mathcal A}$ to a fundamental domain in ${\mathcal T}(S)$ is precompact, 
and the lengths of closed geodesics locally 
separate points in  ${\mathcal T}(S)$. (For example, 
Fenchel-Nielsen  coordinates on  ${\mathcal T}(S)$ 
are determined by lengths of closed geodesics \cite{Buser:1992}.) 
Hence, to prove the theorem we need only discard the possibility 
that there exists a {\em divergent} sequence, 
$[g_n]  \subset {\mathcal A}$, a sequence that 
leaves every compact set in ${\mathcal M}(S)$.

\begin{remark}  \label{InjectivityRemark}
The injectivity radius of $(S, g)$ is determined
by the length spectrum. Indeed, since the curvature is negative, 
the injectivity radius equals half of the length of the shortest 
curve, in other words, half of the infimum of the length spectrum. 
In the case that the number of funnels $n_F=0$, Mumford's lemma
\cite{Mmf71} would give the desired precompactness of (a lift of) 
${\mathcal A}$. If, however, $n_F>0$, then Mumford's lemma
does not apply (directly) to ${\mathcal T}(S)$.
\end{remark}

Recall that a complete hyperbolic surface $(S, g)$ with 
$\chi(S)< \infty$  is determined by its convex core. 
The {\em convex core} can be described
as either the convex hull of the set of closed $g$-geodesics on $S$
or as the hyperbolic surface $(S,g)$ with its funnels removed.

A given cusp end of $(S, g)$ can be compactified by adding one point.  
A closed curve $\alpha \subset S$ will be called 
{\em non-cuspidal} if and only if $\alpha$ is not null-homotopic
after compactifying any single cusp. In other words, 
given a deck group representation  
$\rho: \pi_1(S) \rightarrow SL(2, {\mathbb R})$
for $(S, g)$,  a closed curve $\alpha$
is non-cuspidal if and only if $|{\rm tr}(\rho([\alpha]))| >2$. 
In particular, $\alpha$ is non-cuspidal if and only if 
there exists a (unique) $g$-geodesic that is homotopic to 
$\alpha$.

\begin{lemma}[Geometric Limits]  \label{GeometricLimits}
Let $S$ be a (connected) differentiable surface of 
with $\chi(S')<0$.  Let $[g_n] \in {\mathcal M}(S)$ 
be a divergent sequence with injectivity radius 
uniformly bounded from below by
a positive constant. Then there exist 
\begin{enumerate}
\item a subsequence of metric representatives, $g_n$, 

\item a complete, finitely connected, 
    hyperbolic surface $(R,h)$ with $\chi(R) < \infty$, 

\item a precompact neighborhood $U$ of the convex core of 
      $(R,h)$, and

\item a smooth embedding  $f: U \rightarrow S$, 
\end{enumerate}
such that \renewcommand{\labelenumi}{(\Alph{enumi})}
\begin{enumerate}

\item  each metric $f^*(g_n)$ extends to a complete hyperbolic
       metric $h_n$ on $R$, 

\item  $h_n \rightarrow  h$ in $C^k(U)$ for each 
       $k \in{\mathbb N}$, \label{UConvergence}

\item  for each $n$, the convex core of $(R,h_n)$ lies inside $U$,

\item  given $L>0$, there exists $M>0$ such that if $m>M$
       and $\alpha \subset S$ is a non-cuspidal closed curve that 
       is not homotopic to a curve in $f(U)$, 
       then the $g_m$-length of $\alpha$ is larger than $L$, and 

\item  $\chi(S) < \chi(R)$, \label{EulerInequality}
\end{enumerate}
The surface $(R,g)$ is a {\em geometric limit} of $(S,g_n)$.
\end{lemma}

We postpone the proof of Lemma \ref{GeometricLimits}
to the end of this section.

\vspace{.5cm}

\begin{proof}[Proof of Theorem \ref{IsofiniteTeich}]
By the discussion immediately after the statement of 
Theorem \ref{IsofiniteTeich}, it suffices to show
that there cannot exist a divergent isospectral sequence.
Suppose to the contrary that such a sequence exists.
By Remark \ref{InjectivityRemark}, the injectivity
radius of this sequence is bounded below.

Note that if a connected, complete, hyperbolic surface $(S, g)$ 
satisfies $\chi(S) \geq 0$, then $S$ is either a cylinder 
or the disk. In the latter case, there is only one isometry class. 
In the former, either $S$ is a horn in which case there
is one isometry class, or $S$ is a cylinder in which 
case there is exactly one primitive closed geodesic
whose length---the infimum of the length spectrum---determines 
the isometry class of the metric. In particular, 
any isospectral sequence of metrics on the cylinder 
cannot be divergent.

Thus we assume that $\chi(S)<0$, and hence 
apply  Lemma \ref{GeometricLimits} to obtain 
a geometric limit $(R, h)$. Let $f$, $U$, and $h_n$
be as in Lemma \ref{GeometricLimits}.

Let $\beta$ (resp. $\beta_n$) be the unique $h$-geodesic 
representative (resp. $h_n$-geodesic representative) 
of the homotopy class of a non-cuspidal closed curve 
$\beta \subset R$. We claim that, by passing to a 
further subsequence if  necessary, we may assume that $\beta_n$ 
converges to $\beta$. Indeed, each geodesic 
is a solution to an ordinary
differential equation whose coefficients depend continuously
on the metric. By part (C) of Lemma  \ref{GeometricLimits},
the convex core of each $(R,h_n)$ lies
in the precompact set $U$, and hence there exists a subsequence 
such that the initial conditions converge. Thus, the 
claim follows from part (B) of Lemma \ref{GeometricLimits}
and the continuity of solutions to ordinary
differential equations with respect to coefficients and initial data.

There exists a finite number of homotopy classes of 
simple closed curves,  $\gamma_i$, on $R$ such that 
the $h_n$-lengths (resp. $h$-lengths) of the corresponding 
geodesics determine the metric $h_n$ (resp. $h$) \cite{Buser:1992}.   
By the above, the $h_n$-lengths of the 
$h_n$-representatives of the $\gamma_i$ converge  
to the $h$-length of the $h$-geodesic representing $\gamma_i$.
Since, by assumption, the length spectrum is constant,
the $h_n$-length of each $\gamma_i$ is constant for 
sufficiently large $n$.  Therefore, by passing to a further
subsequence if necessary, we may assume that 
$(R, h_n)$ is isometric to $(R,g)$.

Let $N^0_n(x)$ (resp. $N^{\infty}_n(x)$) denote the number of closed 
$h_n$-geodesics having length less than $x$ that are homotopic 
(resp. not homotopic) to a curve in $f(U)$.   Then we have 
\[    N_{(S, g_n)} (x)~  =~   N^0_n(x) +  N^{\infty}_n(x)~
   =~   N_{(R, h_n)} (x) +  N^{\infty}_n(x). \]
By assumption, the left hand side does not depend on $n$,
and since $(R, h_n)$ is isometric to $(R,h)$ for all $n$,
we have that $N_{(R, h_n)}$ is independent of $n$.
Thus, to obtain a contradiction, it suffices to show 
that  $N^{\infty}_n$ does depend on $n$.

Let $\alpha \subset S$ be a non-cuspidal $g_1$-geodesic
that is not  homotopic to a curve in $f(U)$. We apply
part (D) of Lemma \ref{GeometricLimits} to the curve 
$\alpha$, choosing $L$ to be equal to twice 
the $g_1$-length of $\alpha$. Thus, there exists $M$---necessarily 
larger than $1$---such that the  $g_{M}$-length of any 
non-cuspidal curve that is not homotopic to a curve 
in $f(U)$ is larger than than $L$. In particular, the
$g_M$-length of any closed $g_M$-geodesic 
that is not homotopic to a curve in $f(U)$ is greater
than $L$. 

It follows that for $x$ satisfying $\frac{L}{2} < x < L$,
we have  $N^{\infty}_{M}(x) = 0$ and  $N^{\infty}_{1}(x) >0$. 
Thus, the functions $N^{\infty}_{n}$ depend on $n$, thus  giving
the desired contradiction.
\end{proof}

\vspace{.5cm}

\begin{figure} \label{SetUp}
\caption{The set up for geometric limits}

\psfrag{bK}{$\partial K$}
\psfrag{K2}{$K'$}
\psfrag{gamma}{$\gamma$}
\psfrag{K1}{$K$}
\psfrag{K*}{$K^*$}
\psfrag{U}{$U$}
\psfrag{f(U)}{$f(U)$}
\psfrag{R}{$R$}
\psfrag{D}{$D$}
\begin{center}   \includegraphics{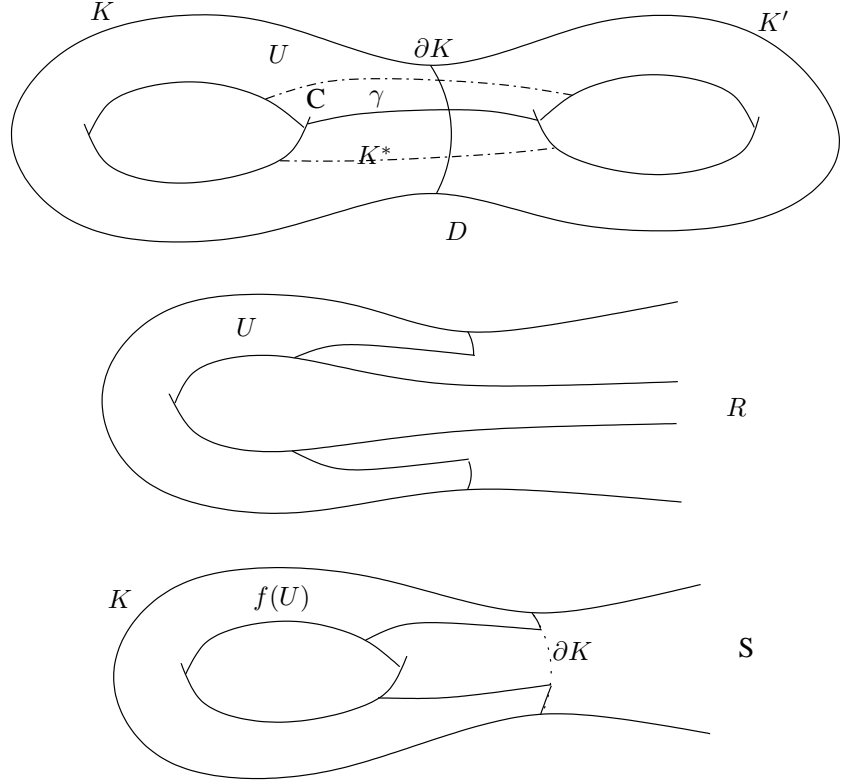} 
\end{center}
\end{figure}

\begin{proof}[Proof of Lemma \ref{GeometricLimits}]
Let $K, K'$ be two isometric copies of the convex core of $(S,g_n)$.
Since the number of funnels and cusps does not depend on $N$,
the diffeomorphism class of these surfaces does not
depend on $n$. Let $D$ be the surface obtained by gluing 
$K$ and $K'$ along their respective boundaries. Since $\partial K$, 
$\partial K'$, are geodesic, the metric $g_n$ extends to
a complete hyperbolic metric on $D$ having finite area. 
Abusing notation slightly, we will refer to this
extended metric as $g_n$.

Note that the divergence of $[g_n]$ implies the divergence
of $[g_n]$ on the double $D$.
Let $\langle g_n \rangle \in {\mathcal T}(D)$
such that $\pi(\langle g_n \rangle) =[g_n]$.
				   
Mumford's lemma \cite{Mmf71} (cofinite case) implies 
that there exists a sequence, $\gamma_n$, of 
multicurves\footnote{Recall that a {\em multicurve} 
is a free homotopy class of a finite disjoint union of 
mutually non-homotopic simple closed curves.
A well-known argument shows that a multicurve has at 
most $3g-3$ components.} on $D$ whose  $\langle g_n \rangle$-length, 
$\ell_n(\gamma_n)$, tends to zero as $n$ tends to infinity.

By Lemma \ref{OrthogonalMulticurve} below, 
each component of $\gamma_n$ meets $\partial K$ orthogonally.  
It follows that the group, $\rm{Mod}_{\partial K}(D)$, 
of mapping classes that preserve (the homotopy class of) $\partial K$
acts on the collection  $\{\gamma_n\}$ with finitely many orbits.
Hence, by taking a subsequence and applying elements
of $\rm{Mod}_{\partial K}(D)$, we may assume that 
$\langle g_n \rangle \subset {\mathcal T}(S)$ was chosen so that
sequence $\gamma_n$ is equal to a constant $\gamma$.

Standard arguments (using, for example, pants decompositions) 
show that there exist representatives $g_n$ of $\langle g_n \rangle$
such that $\partial K$ is $g_n$-geodesic for all $n$,
such that the $g_n$-geodesic representative (still denoted $\gamma$) 
of the multicurve $\gamma$ is independent of $n$, and 
such that $g_n$ converges uniformly on compact subsets
of $D \setminus \gamma$ to a complete (finite area)
hyperbolic metric $g$ on $D \setminus \gamma$.

Let $K^*= K \setminus (K \cap \gamma)$. Then $(K^*,g_n)$ has 
piecewise geodesic boundary with interior angles equal
to $\frac{\pi}{2}$ by Lemma \ref{OrthogonalMulticurve}.
It follows that $K^*$ is geodesically convex and, moreover,
the $g_n$-holonomy (in the sense of hyperbolic structures) 
of any closed curve in $K^*$ cannot be elliptic. Thus,
the hyperbolic surface $(K^*,g_n)$ extends to a complete
hyperbolic surface $(R,h_n)$ with $R$ homeomorphic to 
the interior of $K^*$.

Similarly, we obtain a metric $h$ on $R$ such that $(R,h)$ 
is a complete hyperbolic surface that is an extension of $(K^*,g)$. 
Since the metrics $g_n$ converge to $g$, it follows, 
by analytic continuation, that $h_n$ converges to $h$ 
uniformly on compact subsets of $R$.

Note that since $(K^*,g_n)$ is geodesically convex, 
the convex core of $(R, h_n)$ is bounded by simple closed 
geodesics that lie with $K^*$. 
By Lemma \ref{TransverseCollar} below, there exists,
independent of $n$, a collar neighborhood $C$  of 
$\gamma$  such that no simple closed
$g_n$-geodesic in $K \setminus \gamma$ intersects $C$.
Thus, the convex core of each $(R, h_n)$ lies within
the compact set $U = K^* \setminus C$. Let 
$f: U \rightarrow S$ be the associated inclusion.
Then $R$, $U$, $f$, $h$, and $h_n$ satisfy (A), (B), and (C).
See Figure 2.

By the collar lemma \cite{Buser:1992}, given $L>0$,
there exists $n$ such that if a $g_n$-geodesic on $D$
intersects $\gamma$, then the $g_n$-length must be 
greater than $L$. Part (D) of the Lemma follows.

Let $k$ be the number of components of the multicurve 
$\gamma$. Note that $\partial K^*$ has $4 \cdot k$ corners 
each with  an exterior angle equal to $\pi/2$ by Lemma 
\ref{OrthogonalMulticurve}. Therefore, 
it follows from the Gauss-Bonnet theorem
that 
\begin{equation} \chi(K^*)~ =~ \chi(K)  + k.  \label{GB} \end{equation}
Since  $\chi(R)= \chi(K^*)$ and $\chi(S) = \chi(K)$
and $k>0$, we have part (E).
\end{proof}

\begin{lemma} \label{OrthogonalMulticurve}
There exists a constant $\ell^*$ such that if
$\gamma$ is the $g_n$-geodesic representative of a multicurve
on $D$  of $g_n$-length $\ell(\gamma) < \ell^*$, then each component of $\gamma$ 
intersects  $\partial K$ orthogonally. Moreover, $\ell^*$,
depends only on a lower bound on the injectivity radius of $(D,g_n)$. 
\end{lemma}

\begin{proof}
If  $\ell(\gamma)$  is less than the injectivity radius of $(D, g_n)$,
then no component of $\gamma$ can lie entirely inside one of the two
copies of $K$. Hence each component intersects $\partial K$ 
transversally. To see that this intersection is orthogonal, 
let $\phi: D \rightarrow D$ be the reflection isometry associated 
to the doubling. Then $\phi(\gamma)$ is 
another geodesic multicurve on $D$ that intersects 
$\gamma$ at the fixed point set of $\phi$, $\partial K$.  Since $\phi$ 
is an isometry, $\ell (\phi(\gamma))= \ell (\gamma)$. Recall that 
by the collar lemma, there exists $\ell'$ (independent of the 
complete hyperbolic surface) such that if $\beta_1$ and
$\beta_2$ are intersecting, simple closed geodesics each 
with length less than $\ell'$,
then $\beta_1 =\beta_2$.  

Let $\ell^*$ be the minimum of $\ell'$ and the injectivity radius.
\end{proof}

\begin{figure} \label{PointNoReturnFigure}

\caption{Universal cover of $X$ in Lemma \ref{TransverseCollar}}

\psfrag{a}{$a$}
\psfrag{1/a}{$a^{-1}$}
\psfrag{bS}{$\partial \Sigma$}
\psfrag{bF}{$\partial F$}
\psfrag{gamma}{$\tilde{\gamma}$}
\psfrag{S}{$\Sigma$}
\psfrag{1}{$1$}
\psfrag{el}{$e^{\ell/2}$}
\psfrag{-el}{$-e^{\ell/2}$}
\psfrag{alpha}{$\tilde{\alpha}$}
\psfrag{F}{$F$}
\begin{center}   \includegraphics{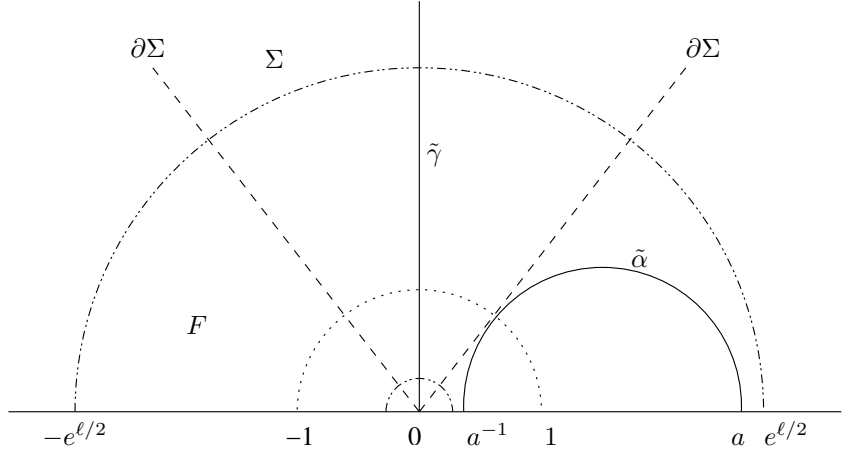} 
\end{center}
\end{figure}

\begin{lemma}[Point of No Return] \label{TransverseCollar}
Let  $\gamma$ be a simple closed geodesic of length $\ell$ 
on a complete hyperbolic surface $X$. If $\alpha$ is 
a simple closed geodesic that does not intersect $\gamma$,
then its distance $d$ from $\gamma$ is at least  
$\log(\coth(\ell/4))$.
\end{lemma}

\begin{proof}
Consider the universal covering of $X$ by the upper half-plane 
$\{(x,y) \in {\mathbb R}^2~ |~ y>0 \}$ such that the
lift of the geodesic $\gamma$ is the $y$-axis: 
$\tilde{\gamma} =\{ (0,y) | y>0 \}$. Let $\tau$ be the 
hyperbolic isometry $\tau(x, y) = (e^{\ell} x, e^{\ell} y)$
stabilizing $\tilde{\gamma}$.  Each tubular neighborhood 
of $\gamma$ is the quotient by $\langle \tau \rangle$ 
of a (Euclidean) sector that is bisected by $\tilde{\gamma}$.

Since both $\gamma$ and $\alpha$ are closed, $d>0$.
Take $\tilde{\alpha}$ to be a lift
of $\alpha$ such that the distance between $\tilde{\gamma}$ 
and $\tilde{\alpha}$ is $d$. (Note that $\tilde{\alpha}$ cannot
be vertical.)

Let $\Sigma$ be the smallest sector bisected by $\tilde{\gamma}$ 
that intersects $\tilde{\alpha}$. (See Figure 3).
Then $\tilde{\alpha}$ meets $\partial \Sigma$ tangentially at one 
point $(x,y)$ and $\Sigma/\langle\tau\rangle$ is a $d$-tubular 
neighborhood of $\gamma$.
Conjugating by a suitable hyperbolic element having axis $\tilde{\gamma}$
and a reflection about $\tilde{\gamma}$ if necessary, 
we may assume without loss of generality, that $x>0$ and $x^2+y^2=1$.
Thus, $\tilde{\alpha}$ is a semi-circle with endpoints 
at $(1/a, 0)$ and $(a,0)$ where $a>1$. 

Consider the following fundamental domain of $\tau$: 
\[ F= \{ (r \cos(\theta), r \sin(\theta) )~  |~    0< \theta < \pi;~  
    e^{-\ell/2} \leq  r \leq e^{\ell/2} \} \]
We claim that $\tilde{\alpha}$ does not intersect either
of the boundary components of $F$: 
$\partial F = \{r=e^{\ell/2}\} \cup \{r= e^{-\ell/2}\}$.  
To see this note that the geodesic $\tilde{\alpha}$ 
meets the (Euclidean) unit circle  $\{r=1\}$ orthogonally 
and hence $\tilde{\alpha}$ is preserved by 
the reflection in $\{r=1\}$. In addition, this 
reflection exchanges the two boundary components of $F$. 
It follows that if $(x',y')  \in \tilde{\alpha} \cap \partial F$, 
then either $\tau(x',y')$ or $\tau^{-1}(x',y')$ belongs 
to $\tilde{\alpha}$.  The corresponding self-intersection 
of $\tilde{\alpha}$ is transverse, a contradiction to the 
simplicity of $\alpha$.

Therefore, we must have 
$e^{-\ell/2} \leq a^{-1} < 1 < a \leq e^{\ell/2}$ 
A straightforward calculation shows that the hyperbolic
width of $\Sigma$ is at least $2 \log(\coth(\ell/4))$. 
\end{proof}



\appendix
\section{Special function calculations}\label{sp.fun.sec}
To calculate $\Phi_F(s)$ and $\Phi_H(s)$ explicitly, 
we derive and use an explicit representation 
of the integral kernel $G_{\bbH}(s;z,z')$  of the resolvent
operator $(\Delta _{{\bbH}}-s(1-s))^{-1}$.  In particular,  if we denote
\begin{equation*}
     \sigma(z,z')~  =~  \cosh^2(d(z,z')/2) =\frac{(x-x')^{2}+(y+y')^{2}}{4yy'},
\end{equation*}
then $G_{\bbH}(s,z,z') = g_s(\sigma)$ where $g_s$ can be expressed 
explicitly in terms of a hypergeometric function: 
\begin{equation*}
g_s(\sigma)  = \frac{1}{4\pi} 
\frac{\Gamma(s)^2}{\Gamma (2s)}\sigma ^{-s}F(s,s;2s,\sigma ^{-1}),
\end{equation*}
where 
\begin{equation*}
F(a,b;c;z)=1+\frac{ab}{c}z+\frac{a(a+1)b(b+1)}{2!c(c+1)}z^{2}+\dots 
\end{equation*}
(see, for example, formula (6.6.1) of \cite{He:1983}).
Euler's representation of the hypergeometric function,
$$
F(a,b;c;z) = \frac{\Gamma(c)}{\Gamma(b) \Gamma(c-b)} 
\int_0^1 \frac{t^{b-1} (1-t)^{c-b-1}}{(1-tz)^a} \>dt,
$$
gives us the useful representation:
\begin{equation}\label{gs.formula}
G_{\bbH}(s,z,z') = g_s(\sigma) = 
\frac1{4\pi} \int_0^1 \frac{t^{s-1} (1-t)^{s-1}} {(\sigma-t)^s} \>dt,
\end{equation}
for $\sigma>1$, $\re(s) > 0$.

As an application, we can evaluate
$$
\lim_{\sigma\to1} \frac{d}{ds} g_s(\sigma) = \frac{1}{4\pi} \int_0^1 t^{s-1} 
\log(t(1-t))\>dt
= -\frac{1}{2\pi} \Psi'(s).
$$
Since $(d/ds) R_\bbH(s) = -(2s-1) R_\bbH(s)^2$ this implies that
\begin{equation}\label{rsq.diag}
(R_\bbH(s)^2)(s,w,w) = \frac{1}{2\pi} (2s-1)^{-1} \Psi'(s)
\end{equation}
By combining this with the same expression evaluated at $1-s$ and using the 
identity $\Gamma(s)\Gamma(1-s) = \pi \csc (\pi s)$, we can also derive
$$
\frac{d}{ds} [G_\bbH(s,z,w) - G_\bbH(1-s,z,w)]_{z=w} = -\frac{\pi}{2} \csc^2 (\pi s).
$$
This can be integrated to give 
\begin{equation}\label{ghgh.diag}
[G_\bbH(s,z,w) - G_\bbH(1-s,z,w)]_{z=w} = \frac12 \cot(\pi s)
\end{equation}

\subsection{Funnel}
As above, $F$ is the model funnel (hyperbolic half-cylinder) and
$\Delta _{F}$ denotes the Laplacian with Dirichlet boundary
conditions.  $F$ is half of the cylinder $M = \Gamma_\ell\backslash\bbH$.
 The resolvent kernel for $\Delta _{F}$ is obtained from that of $\Delta _{M}$ by
the method of images: 
\begin{equation}\label{res.images}
       G_{F}(s;z,w)=G_{M}(s;z,w)-G_{M}(s;z,\tau w),
\end{equation}
where $\tau$ denotes the involution $(x,y)\mapsto (-x,y)$ and this is really 
a statement about the lifts of the kernels to $\bbH\times \bbH$.   

\begin{proof}[Proof of Proposition \ref{funnel.phi}]
Using (\ref{res.images}) the integrand in the definition on $\Phi_F(s)$  
can be written 
\begin{equation}
\begin{split}
\lbrack G_{F}(s;z,w)-G_{\bbH}(s;z,w)]_{w=z}& =\sum_{m\neq0}
g_{s}\left( \frac{(x-e^{m\ell}x)^{2}+(y+e^{m\ell}y)^{2}}
{4e^{m\ell}y^{2}}\right)  \\
& \quad -\sum_{m\in \mathbb{Z}}g_{s}\left( \frac{(x+e^{m\ell}x)^{2}+
(y+e^{m\ell}y)^{2}}{4e^{m\ell}y^{2}}\right) 
\end{split}
\label{eq.integrand}
\end{equation}

With an inverse Mellin transform, Patterson \cite{Patterson:1989} 
computes that 
\begin{equation}\label{gs.int}
(2s-1)\int_{-\infty}^\infty g_{s}\left( (1+t^{2})
\frac{(1+\kappa )^{2}}{4\kappa }
\right) \>dt=\frac{\max (\kappa ,\kappa ^{-1})^{1/2-s}}{\kappa ^{1/2}+
\kappa^{-1/2}}.
\end{equation}
We will apply this by choosing fundamental domain $\mathcal{F} =
\bbR_+ \times [1,e^\ell]$.
When (\ref{eq.integrand}) is integrated over $\mathcal{F}_\lambda$, 
the $m=0$ term in the second sum contributes 
\begin{equation*}
\begin{split}
-(2s-1)\int_{1}^{e^\ell}\int_{0}^{\infty }g_{s}(1+\tfrac{x^{2}}{y^{2}})
\frac{dx\>dy}{y^{2}}& =-\frac{2s-1}{2}\int_{1}^{e^\ell}\int_{-\infty}^\infty
g_{s}(1+t^{2})\frac{dt\>dy}{y} \\
& =-\frac{\ell}{4}.
\end{split}
\end{equation*}
For the $m\neq 0$ the appropriate substitution is 
\begin{equation*}
t=\frac{(1\mp e^{m\ell})}{(1+e^{m\ell})}\frac{x}{y},
\end{equation*}
(with $\mp $ depending on whether the involution $\tau$ is inserted or not).  
Integrating an $m\ne 0$ from (\ref{eq.integrand}) over $\mathcal{F}_\lambda$ 
then gives
\begin{equation*}
\begin{split}
&(2s-1)\int_{\mathcal{F}}g_{s}\left( \frac{(x\mp e^{m\ell}x)^{2}
+(y+e^{m\ell}y)^{2}}{4e^{m\ell}y^{2}}\right) \>\frac{dx\>dy}{
y^{2}} \\
& \quad =\frac{\ell}{2} \cdot \frac{(1+e^{m\ell})}{(1\mp e^{m\ell})}
(2s-1)\int_{-\infty}^\infty g_{s}\left( (1+t^{2})\frac{(1+\lambda ^{2})^{2}}{4
e^{m\ell}}\right) \>dt \\
& \quad =\frac{\ell}{2} \cdot \frac{(1+e^{m\ell})}{(1\mp e^{m\ell})}
\frac{e^{(1/2-s)|m|\ell}}{e^{m\ell/2}+e^{-m\ell/2}} \\
& \quad =\frac{\ell}{2} \cdot \frac{e^{-s|m|\ell}}{1\mp e^{-|m|\ell}}
\end{split}
\end{equation*}
Thus the full contribution from the $m\neq 0$ terms is 
\begin{equation*}
\begin{split}
\frac{\ell}{2} \sum_{m\neq 0}\left( \frac{e^{-s|m|\ell}}{
1-e^{-|m|\ell}}-\frac{e^{-s|m|\ell}}{1+e^{-|m|\ell}}\right) &
=2\ell \sum_{m\geq 1}\left( \frac{e^{-(s+1)m\ell}}{1-
e^{-2m\ell}}\right) \\
& =2\frac{d}{ds}\sum_{k\geq 1}\log (1-e^{-(s+2k+1)\ell})
\end{split}
\end{equation*}
Thus, taking into account all of the terms in (\ref{eq.integrand})
gives
$$
(2s-1)\int_{\mathcal{F}}[G_{F}(s;z,w)-G_{\mathbb{H}
^{2}}(s;z,w)]_{w=z}\>dg(z)
 =  -\frac{\ell}{4} + \frac{d}{ds}\sum_{k\geq 1}\log (1-e^{-(s+2k+1)\ell})^2
$$
\end{proof}


\subsection{Horn}

The resolvent and generalized eigenfunctions for the horn $H = 
\Gamma_\infty\backslash  \bbH$ may be
obtained by summing $G_\bbH$ over $\Gamma_\infty$.   This makes it 
possible to compute $\Phi_H(s)$ fairly directly.

\begin{proof}[Proof of Proposition \ref{horn.phi}]
For the proof it suffices to assume $\re(s)>1/2$.
The integrand in (\ref{eq.horn.trace}) can be written 
\begin{eqnarray*}
[G_{H}(s;z,w)-G_{\mathbb{H}}(s;z,w)]_{w=z}
&=&  \sum_{k\ne 0} G_{\mathbb{H}}(s,z,z+k) \\
&=&  2 \sum_{k=1}^\infty g_s\left(1 + \frac{k^2}{4y^2}\right).
\end{eqnarray*}
We will use $[0,1]\times \bbR_+$ as our fundamental domain for $H$.
The integral in the definition of $\Phi_H(s)$ is well-behaved at the funnel end, but
needs regularization at the cusp end:
$$
\Phi_H(s) 
= 2 (2s-1) \FPe \int_0^{1/\varepsilon}\int_0^1 \sum_{k=1}^\infty 
g_s\left(1 + \frac{k^2}{4y^2}\right) \frac{dx\>dy}{y^2}.
$$
Substituting $u = \frac{k}{2y}$, this becomes
\begin{eqnarray*}
\Phi_H(s) 
&=& 4 (2s-1) \FPe \sum_{k=1}^\infty \int_{k\varepsilon/2}^\infty  
g_s(1 + u^2) \frac{1}{k}\>du \\
&=& 4 (2s-1) \FPe \int_{0}^\infty  g_s(1 + u^2) \left[\sum_{1\le k 
\le 2u/\varepsilon} \frac{1}{k}\right] \>du.
\end{eqnarray*}
In order to extract the finite part, we use
$$
\sum_{1\le k \le 2u/\varepsilon} \frac{1}{k} = \log \frac{2u}{\varepsilon} + \gamma
+ O(\sqrt{\varepsilon/2u}),
$$
where $\gamma$ is Euler's constant.  The $O(\sqrt{\varepsilon/2u})$ contribution
will vanish in the $\varepsilon\to0$ limit, and $\log\varepsilon$ is dropped in 
taking the finite part, leaving
\begin{equation}\label{phih.eq}
\begin{split}
\Phi_H(s) &=  (\log 2 + \gamma) \cdot 4(2s-1) \int_0^\infty g_s(1+u^2) \>du \\
&\quad+ 4(2s-1) \int_0^\infty g_s(1+u^2) \log u\>du.  
\end{split}
\end{equation}
The formula (\ref{gs.int}) gives:
\begin{equation}\label{gs.int2}
4(2s-1) \int_0^\infty g_s(1+u^2) \>du = 1.
\end{equation}
For the other integral in (\ref{phih.eq}) we use (\ref{gs.formula}) to write
$$
4(2s-1) \int_0^\infty g_s(1+u^2) \log u\>du 
= \frac{2s-1}{\pi} \int_0^\infty \int_0^1 \frac{t^{s-1} (1-t)^{s-1}}{(1-t+u^2)^s}
\log u\>du\>dt.
$$
Substituting $u = \sqrt{(1-t)(1-w)/w}$ yields
\begin{equation*}
\begin{split}
&4(2s-1) \int_0^1 g_s(1+u^2) \log u\>du \\
&\quad= \frac{2s-1}{4\pi} \int_0^1 \int_0^1 t^{s-1} 
(1-t)^{-1/2} w^{s-3/2} (1-w)^{-1/2} \log [(1-t)(1-w)/w]\>dw\>dt.
\end{split}
\end{equation*}
These are straightforward Beta-type integrals, whose evaluation gives:
$$
4(2s-1) \int_0^1 g_s(1+u^2) \log u\>du  
= -\gamma - \log 4 - \Psi(s+1/2) + \frac{1}{2s-1}.
$$
Using this last equation together with (\ref{gs.int2}) in (\ref{phih.eq}) gives
the stated result.
\end{proof}

The preceding calculation shows that 
$\res_{s=1/2} \Phi_H(s) = 1/2$.  This residue can be
identified as originating in the funnel end of $H$, with no 
contribution from the cusp.
To make this statement precise, we introduce a cutoff $\eta\in\cinf(\bbR_+)$
which is equal to 0 in a neighborhood of 0 and equal to 1 in a neighborhood of 
$\infty$.  Note that the integrand used to define $\Phi_H(s)$ is explicitly given by
$$
\varphi_H(s;y) = 2(2s-1) \sum_{k=1}^\infty g_s\left(1+ \frac{k^2}{4y^2}\right).
$$
\begin{lemma}\label{cusp.res}
The contribution from the cusp end to the $0$-trace in $\res_{s=1/2} \Phi_H(s) = 1/2$ 
is zero, in the sense that:
$$
\res_{s=1/2} \left[ \FP_{\varepsilon\to0} \int_0^{1/\varepsilon} \eta(y) 
\varphi_H(s;y) \frac{dy}{y^2}\right] = 0.
$$
\end{lemma}
\begin{proof}
It is equivalent to show that the $0$-trace comes exclusively from the funnel end,
i.e. that
$$
\res_{s=1/2} \left[ \FP_{\varepsilon\to0} \int_{\varepsilon}^\infty (1-\eta(y)) 
\varphi_H(s;y) \frac{dy}{y^2}\right] = \frac12.
$$
Using the hypergeometric series representation of $g_s(\sigma)$, we can deduce 
that $(1-\eta(y)) \varphi_H(s;y)$ is of the form $y^{2s} f_s(y)$, with $f_s \in
\cinf_c(\bbR_+)$ and
$$
f_s(0) = \frac{(2s-1)}{2\pi} \frac{\Gamma(s)^2}{\Gamma(2s)} 4^s \zeta(2s).
$$
Integration by parts shows that, for $s\ne 1/2$,
$$
\FP_{\varepsilon\to0} \int_{\varepsilon}^\infty  
y^{2s-2} f_s(y) \>dy = -\frac1{(2s-1)} \int_0^\infty y^{2s-1} f'_s(y)\>dy.
$$
Thus 
$$
\res_{s=1/2} \left[ \FP_{\varepsilon\to0} \int_0^{1/\varepsilon} (1-\eta(y)) 
\varphi_H(s;y) \frac{dy}{y^2}\right] = \frac12 f_{1/2}(0)
$$
The result then follows from $f_{1/2}(0) = 1$, which holds because
the factor $(2s-1)$ exactly cancels the pole of $\zeta(2s)$ at $s=1/2$.
\end{proof}




\end{document}